\documentclass[12pt]{article}
\scrollmode

\usepackage{amsfonts}
\usepackage{amsmath}
\usepackage{amssymb}
\usepackage{graphicx}

\allowdisplaybreaks

\title{K-processes, scaling limit and aging for the REM-like trap model
}
\date{}
\author{Luiz Renato Fontes \footnote{IME-USP, Rua do Mat\~ao 1010, 05508-090
S\~ao Paulo SP,  Brazil, lrenato@ime.usp.br}
\and
Pierre Mathieu \footnote{CMI, 39 rue Joliot-Curie, 13013 Marseille, France,
pierre.mathieu@cmi.univ-mrs.fr}      }

\def\eps{\varepsilon}
\def\qed{\hfill\rule{.2cm}{.2cm}}
\def\P{{\mathbb P}}
\def\esp{{\mathbb E}}
\def\Z{{\mathbb Z}}
\def\N{{\mathbb N}^\ast}
\def\bN{\bar{\mathbb N}^\ast}
\def\R{{\mathbb R}}

\def\1{{\mathbf 1}}

\def\a{\alpha}
\def\o{\omega}  
\def\l{\lambda}
\def\La{\Lambda}
\def\tl{\tilde\lambda}
\def\tln{\tilde\lambda_n}
\def\g{\gamma}
\def\gx{\Gamma^{(x)}}
\def\gu{\Gamma^{(1)}}
\def\gi{\Gamma^{(i)}}
\def\gz{\Gamma^{(0)}}
\def\x{\Gamma}
\def\G{\Gamma}
\def\s{\sigma}
\def\sx{\s^{(x)}}
\def\e{\epsilon}
\def\d{\delta}

\def\tix{T^{(x)}}

\def\n{{\cal N}}

\def\Ga{{\cal G}}

\def\E{{\cal E}}

\def\A{{\cal A}}

\def\T{{\cal T}}
\def\RR{{\cal R}}
\def\S{{\cal S}}
\def\C{{\cal C}}
\def\D{{\cal D}}

\def\x1{X_1^{(1)}}

\def\txc{\tilde X^c}
\def\txcy{\tilde X^{c,y}}
\def\txci{\tilde X^{c,\infty}}
\def\tx0{\tilde X^0}
\def\tmt{\theta_{m,t}}

\def\nn{\nonumber}

\def\={&=&}
\def\+{&+&}

\newtheorem{theo}{Theorem}[section]
\newtheorem{prop}[theo]{Proposition}
\newtheorem{lm}[theo]{Lemma}
\newtheorem{cor}[theo]{Corollary}
\newtheorem{rmk}[theo]{Remark}
\newtheorem{df}[theo]{Definition}
\def\beq{\begin{equation}}
\def\eeq{\end{equation}}
\newcommand{\bei}{\begin{itemize}}
\newcommand{\eei}{\end{itemize}}
\newcommand{\ben}{\begin{enumerate}}
\newcommand{\een}{\end{enumerate}}
\newcommand{\beqn}{\begin{eqnarray}}
\newcommand{\beqnn}{\begin{eqnarray*}}
\newcommand{\eeqn}{\end{eqnarray}}
\newcommand{\eeqnn}{\end{eqnarray*}}
\newcommand{\brm}{\begin{rmk}}
\newcommand{\erm}{\end{rmk}}

\begin{document}\maketitle


\begin{abstract}
We study K-processes, which are Markov processes in a denumerable
state space, all of whose elements are stable, with the exception
of a single state, starting from which the process
enters finite sets of stable states with uniform distribution. We
show how these processes arise, in a particular instance, as
scaling limits of the REM-like trap model
``at low temperature'', and subsequently derive aging results for
those models in this context.

\end{abstract}


\smallskip




\section{Introduction}

In this paper, we study some properties of a family of Markov
processes, which we call K-processes, in particular, and that's
our main motivation, its relationship in a special case with the
scaling limit of a trap model associated to the Random
Energy Model (REM) at low temperature, as well as with the aging
phenomenon exhibited by that model~\cite{kn:BD}. These processes
are thus prototypes of infinite volume dynamics for low
temperature (mean field) spin-glasses.

They have the following remarkable characteristic property. Their
state space is denumerably infinite (we take it to be
$\{1,2,\ldots,\infty\}$), with a single {\em unstable} state, where
by unstable we mean that the process spends $0$ time at that
state at each visit to it; as we'll see, that state may be
either {\em instantaneous} or {\em fictitious} (which are standard
terms) in different cases.
When in a stable state, the process waits for an exponential time
and jumps to the unstable state, starting from which, and
here's the striking feature, it enters any finite set of stable
states with uniform distribution. In the context of spin glasses,
the stable states represent the low energy configurations, and the
unstable state represents the high energy configurations.
The apparent paradox of the uniformity property is elucidated by
a summability condition on the inverse of the jump rates.

It turns out that a class of processes with this uniformity
property was introduced by Kolmogorov as an example of a Markov
process with an instantaneous state, thus not satisfying his
equations in their usual form~\cite{kn:K}, and it has subsequently
been considered by many authors. This class comprises all members
of the family we study in the present paper but for an important
special case, precisely the one related to the trap model.
See Remark~\ref{rmk:kp} below for more details.

We have two approaches: an analytical one, based on Dirichlet
forms, introduced in Section~\ref{sec:anal}; and one based on an
explicit probabilistic construction, in Section~\ref{sec:prob}, at
the end of which we argue the equivalence of both points of view.
In Section~\ref{sec:char}, we derive a characterization
result for K-processes. Section~\ref{sec:bou} is devoted to the
scaling limit of the REM-like trap model, and to deriving an aging
result for the associated K-process in this context, which can be
seen as an aging result for the trap model itself.

Aging is a dynamical phenomenon observed in disordered systems
like spin-glasses at low temperature, signaled by the existence
of a limit of a given 
two-time correlation function of the system started at a high
temperature configuration/state, as both times diverge keeping a
fixed ratio between them; the limit should be a nontrivial
function of the ratio. This is thus a far-from-equilibrium phenomenon.
It has been observed in real spin glasses and studied extensively in
the physics literature. See~\cite{kn:BCKM} and references therein.

In~\cite{kn:BD}, a phenomenological model for a Glauber dynamics
for the Random Energy Model (REM) is introduced, namely the trap
model (in the complete graph), and an aging result for that model
established. See more on the the trap model and what is meant by
an aging result in Section~\ref{sec:bou}. Roughly speaking, the
trap model is a symmetric continuous time random walk, typically
in a regular graph, finite or infinite. The jump rates at the
vertices are i.i.d.~random variables with a polynomial tail at the
origin, whose degree is related to temperature, so that
degree less than 1 is equivalent to low temperature. We'll assume
this regime throughout.

In the mathematics literature, much attention has recently been
given to the trap model, and many aging results were derived for it.
In~\cite{kn:ABG1,kn:ABG2}, the trap model in the hypercube is
studied, with the rates given by energies of the REM associated to
the vertices of the hypercube. The aging result obtained
in~\cite{kn:ABG2} is for the same correlation function as one
considered in~\cite{kn:BD} with the same limit, in this fashion
giving support to the phenomenology underlying the adoption of the
trap model by the authors of the latter paper. In~\cite{kn:BF},
among other results, the aging result in~\cite{kn:BD} aluded to above
is established in a mathematically rigorous fashion.

The trap model in $\Z$ was considered in~\cite{kn:FIN}
and~\cite{kn:AC}; the one in $\Z^2$, in~\cite{kn:ACM,kn:C}; the
one in $\Z^d$, $d\geq3$ in~\cite{kn:C}. In~\cite{kn:AC1} a
comprehensive approach to obtaining aging results for the
trap model on a class of graphs, including
$\Z^d$ and tori in two and higher dimensions, the complete graph,
the hypercube, is developed.

In most of the above cited work, aging is derived for given
correlation functions, without specific regard to the fact that
aging may arise as a scaling property of the full dynamics. As
in~\cite{kn:FIN} and~\cite{kn:AC}, we follow the latter approach
for the trap model in the complete graph, and derive its scaling
limit (see Theorem~\ref{teo:scabou} below); aging results follow
(after a further limit is taken, as explained below;
see also Theorem~\ref{teo:age} and Corollary~\ref{cor:age} below).

It should be noted that, since a time divergence is involved, the
scaling limit of the rates (or alternatively the average holding times)
should be taken together with the scaling limit of the dynamics,
the limiting object acting as a disordered set of parameters for
the limiting dynamics. The rescaling is of time only (in such a
way that the lowest rates are of order 1), since space isn't
relevant for the model in the complete graph. The scaling limit
results as roughly speaking a dynamics in the deepest traps (but
the remainder states play a role: they're lumped together in the
limit in a single unstable state).

In this model, in the scaling limit, aging is a phenomenon of the
dynamics at {\em vanishing} times: at order 1 or larger times the
dynamics is close enough to or in equilibrium, in contrast to the
one dimensional case of~\cite{kn:FIN} and~\cite{kn:AC}, where it
could be said that aging occurs for fixed macroscopic times. This
should be compared to the aging result in~\cite{kn:BD}
and~\cite{kn:BF}, alluded to above, which takes place in a large
microscopic time regime (in our case, it occurs at short
macroscopic times), and also to the aging result of~\cite{kn:AC1}
for the complete graph, taking place at {\em mesoscopic} time
scales; as far as the three regimes can be compared, they
coincide, perhaps not surprisingly. See Remarks~\ref{rmk:tto0} 
and~\ref{rmk:ac}.

By taking the scaling limit first, and the aging limit
after, we can see aging as a macroscopic phenomenon (taking place
in the limiting dynamics). We point out that the latter limit holds
for almost every realization of the underlying (macroscopic) disorder:
Theorem~\ref{teo:age} and Corollary~\ref{cor:age} are almost sure
aging results.

The scaling limit for the trap model isn't relevant only as a
background for aging, even though that's our main motivation for
taking it in this paper. It contains also information about other
important features of the dynamics at long microscopic times: from
aging at short macroscopic times, to approach to equilibrium at
large macroscopic times. So it has an interest of its own.
Inasmuch as the REM is a prototype for a (mean-field) spin glass,
and the trap model in the complete graph is a prototype for a
Glauber dynamics for the REM at low temperature, this scaling
limit turns up as a prototype for an infinite-volume dynamics of a
(mean-field) spin glass at low temperature. We expect the same
process to arise as an appropriate scaling limit for the trap
model in the hypercube (as dimension diverges), and also for the
hopping dynamics for the REM, either in the complete graph or the
hypercube. It is conceivable that it will also be the scaling
limit of the Metropolis dynamics for the REM in the hypercube
(see, e.g.,~\cite{kn:FIKP} for a definition of this dynamics). We
also expect variants of the K-process to show up as scaling limits
for dynamics of other mean-field models at low temperature, like
the GREM, and that they will also exhibit aging.

Our first step in this study is to describe the class of processes
that arise as the scaling limit of the trap model in the complete
graph. Since they are closely related to
the above mentioned class of processes introduced by Kolmogorov
through the above mentioned uniformity property, which turns out
to characterize the family consisting of both classes (see
Section~\ref{sec:char} and Theorem~\ref{teo:char}), we chose to
start by defining, constructing and studying relevant properties
of that larger family, which we refer to as K-processes.

As mentioned above, we do that analitically, through the Dirichlet
form associated to the process (in Section~\ref{sec:anal}), and,
alternatively, through a probabilistic construction (in
Section~\ref{sec:prob}). The former way has the advantage that the
K-processes (are reversible and) have quite simple Dirichlet
forms, which facilitate the analysis of quantities like the Green
function (see Subsections~\ref{ssec:hit} and~\ref{ssec:ext}).

The probabilistic construction, besides having its own interest,
allows for a direct analysis of the scaling limit for the trap
model and the aging issue, without the need of taking
transforms\footnote{But we do rely on a Tauberian theorem at a
specific point of our argument; see proof of
Theorem~\ref{teo:age}.}, and entails the inclusion of more general
aging functions in the analysis and results (see
Theorem~\ref{teo:age} and Corollary~\ref{cor:age}), at little
extra effort. See Section~\ref{sec:bou} and
Subsection~\ref{ssec:age}.

The analytical construction also leads to simple derivations of aging
results in a weak sense, after taking Laplace transforms. See beginning
of the proof of Theorem~\ref{teo:age}.

In connection with another area of research, as we briefly discuss
in Remark~\ref{fuku} in Section~\ref{sec:char}, a K-process
can be viewed as a one-point extension of a Markov process beyond
its killing time, an object which is of current
interest~\cite{kn:FT,kn:CFY}.



\section{Dirichlet forms approach}
\label{sec:anal}

\subsection{Construction}

Let $\bN$ the be one point compactification of $\N=\{1,2,\ldots\}$, with $\infty$ denoting
the extra point. In other words, we take $\bN$ with any fixed metric $d$ making it compact.
For definiteness, take
\begin{equation}\label{eq:met}
    d(x,y)=\left|x^{-1}-y^{-1}\right|,\,x,y\in\bN
\end{equation}
(with $\infty^{-1}=0$).

Let $\g:\N\to(0,\infty)$ be such that
\begin{equation}\label{eq:sum}
  \sum_{x\in\N}\g(x)<\infty.
\end{equation}
We extend $\g$ to $\bN$ by declaring
\begin{equation}\label{eq:ginf}
    \g(\infty)=0.
\end{equation}

Let $\C$ be the space on continuous real valued functions on $\bN$ and define
\beqn
\D=\{ f:\bN\rightarrow \R\; s.t.\, \sum_x (f(x)-f(\infty))^2<\infty\}\,.
\eeqn
('$\sum_x$' usually stands for '$\sum_{x\in\N}$'.)
Note that $\D$ is a dense subset of $\C$.

For $f,g\in\D$, consider the bilinear symmetric form
\beqn
\E(f,g)=\sum_x (f(x)-f(\infty))(g(x)-g(\infty))\,.
\eeqn

\begin{lm}
\label{lm:diri} $(\E,\D)$ is a regular Dirichlet form acting on
$L^2(\bN,\g)$ in the sense of~\cite{kn:FOT}.
\end{lm}

\noindent {\bf Proof} First note that $\g$ has {\em full support} since we have assumed that
$\gamma(x)>0$ for all $x\in\N$.

Clearly $\E$ is bilinear and symmetric. We should check that
$\D\subset L^2(\bN,\g)$: let $f\in\D$. W.l.o.g.~assume that
$f(\infty)=0$. Therefore $\E(f,f)=\sum_x f(x)^2<\infty$ and
$\sum_x f(x)^2\g(x)\leq (\sup_x\g(x))\sum_x f(x)^2<\infty$.

It is easy to check that contractions act on $\E$ so that $\E$ is a Markovian form.

The last point is to prove that $\D$ is {\em complete} for the norm induced by the bilinear form
$\E$: assume that $f_n\in\D$ satisfies $f_n\rightarrow 0$ in $L^2(\bN,\g)$ and
$\E(f_n-f_m,f_n-f_m)\rightarrow 0$ as $n$ and $m$ tend to $\infty$. Then we must have
$f_n(x)\rightarrow 0$ for any $x\in\N$ (because $\g(x)>0$).
Also, for any $\eps>0$ there exists $n_0$ s.t. for any $n,m\geq n_0$ and any $x\in\N$,
\beqnn \vert f_n(x)-f_n(\infty)-f_m(x)+f_m(\infty)\vert\leq \eps\,.\eeqnn
(This comes from the assumption $\E(f_n-f_m,f_n-f_m)\rightarrow 0$.)
Letting $m$ go to $\infty$ and then $x$ go to $\infty$, we get that
$\limsup_m \vert f_m(\infty)\vert\leq\eps$ and therefore $f_m(\infty)\rightarrow 0$ as $m$ tends to $\infty$.
By Fatou's Lemma,
\beqnn
&&\E(f_n,f_n)=\sum_x (f_n(x)-f_n(\infty)^2=\sum_x\liminf_m (f_n(x)-f_m(x)-f_n(\infty)+f_m(\infty))^2\\
&\leq& \liminf_m \sum_x(f_n(x)-f_m(x)-f_n(\infty)+f_m(\infty))^2=\liminf_m \E(f_n-f_m,f_n-f_m)\,,
\eeqnn
and therefore $\E(f_n,f_n)\rightarrow 0$. \qed
\begin{rmk}
\label{rmk:mark} We conclude from the above lemma that there
exists a strong reversible Markov process, in fact a Hunt process,
whose Dirichlet form is $(\E,\D)$ on $L^2(\bN,\g)$. We shall
denote it by $(\P_x,x\in\bN),(P_t,t\geq 0)$.
\end{rmk}

\subsection{Computation of hitting times and capacities}
\label{ssec:hit}

Given the explicit enough form of $\E$ it is easy to compute the law of some hitting times and entrance
laws.
\begin{lm}
\label{lm:exi}
Let $\tau_x=\inf\{t\, ;\, X(t)\not=x\}$.
Then
\begin{equation}
\label{eq:exi}
\esp_x(e^{-\lambda\tau_x})=\frac 1{1+\lambda\g(x)}\,.
\end{equation}
\end{lm}

\noindent {\bf Proof} The function $y\rightarrow \esp_y(e^{-\lambda\tau_x})$ is the
minimizer of the expression
$\E(u,u)+\lambda \g(u^2)$ among functions $u$ satisfying $u(y)=1$ for $y\not=x$.
But, for such a function $u$, we have
$\E(u,u)+\lambda \g(u^2)=(u(x)-1)^2+\lambda \g(x)u(x)^2+\lambda (1-\g(x))$ that is minimal
for $u(x)=\frac 1{1+\lambda\g(x)}$. \qed

\begin{lm}
\label{lm:ent}
Let $\s_\infty=\inf\{ t\, ;\, X(t)=\infty\}$. Then
\beqn \esp_x(e^{-\lambda\sigma_\infty})=\frac 1{1+\lambda\g(x)}\,.\eeqn
\end{lm}

\noindent {\bf Proof} We now have to minimize
$\E(u,u)+\lambda \g(u^2)$ among functions $u$ satisfying $u(\infty)=1$.
But for such a function $u$, we have
$\E(u,u)+\lambda \g(u^2)=\sum_x (u(x)-1)^2+\lambda\g(x)u(x)^2$ that is minimal
for $u(x)=\frac 1{1+\lambda\g(x)}$. \qed

\begin{rmk}
\label{rmk:mark1}
In particular note that $\s_\infty<\infty$ $\P_x$.a.s. Hence $\P_\infty$ is well defined.
Since $\esp_x(e^{-\lambda\sigma_\infty})=\esp_x(e^{-\lambda\tau_x})$ and
$\tau_x\leq\sigma_\infty$, we
must have $\tau_x=\sigma_\infty$ $\P_x$.a.s. In particular $X(\tau_x)=\infty$ $\P_x$.a.s.
\end{rmk}

\begin{lm}
\label{lm:enta} Let $A$ be a finite subset of $\N$ of size $n$,
and $\tau^{A}=\inf\{ t\, ;\, X(t)\in A\}$. Then, for any function
$f: A\rightarrow \R$, any $\lambda>0$ and any $y\notin A$, we have
\beqn \esp_y(f(X_{\tau^{A}})e^{-\lambda\tau^{A}}))
=\frac1{n+(1+\lambda\g(y))\sum_{x\notin
A}\frac{\lambda\g(x)}{1+\lambda\g(x)}}\,
 \sum_{x\in A} f(x)
\,.\eeqn
In particular, for $\lambda=0$, we find that the law of $X(\tau^{A})$ is uniform over
$A$.
\end{lm}

\noindent {\bf Proof} We have to minimize
$\E(u,u)+\lambda \g(u^2)$ among functions $u$ satisfying
$u(x)=f(x)$ for $x\in A$. For such a function
$\E(u,u)+\lambda \g(u^2)=
\sum_{x\in A} (f(x)-u(\infty))^2+\sum_{x\notin A} (u(x)-u(\infty))^2+ \lambda\g(x)u(x)^2
+\lambda\sum_{x\in A}\g(x)f(x)^2$. The solution has the form
$u(y)=\frac {u(\infty)}{1+\lambda\g(y)}$ for $y\notin A$ and we find $u(\infty)$ by minimizing
$\sum_{x\in A} (f(x)-u(\infty))^2+ u(\infty)^2
\sum_{x\notin A}\frac{\lambda\g(x)}{1+\lambda\g(x)}$. \qed

\medskip

After a similar computation, we get the following.
\begin{lm}
\label{lm:enta1} Let $A$ be as in the previous lemma. Then
\beqn \esp_\infty(f(X_{\tau^{A}})e^{-\lambda\tau^{A}}))
=\frac1{ n+\sum_{x\notin A}\frac{\lambda\g(x)}{1+\lambda\g(x)}}\,\sum_{x\in A} f(x)
\,.\eeqn
\end{lm}

It is also possible to compute the Green kernel
\begin{equation}
  \label{eq:greenk}
  g_\lambda(x)=\lambda\int_0^\infty e^{-\lambda s}\,\P_\infty(X(s)=x)\,ds.
\end{equation}
 The Markov property gives:
\begin{equation}
  \label{eq:marpro}
g_\lambda(x)=\esp_\infty(e^{-\lambda\tau^{\{x\}}})(1-\esp_x(e^{-\lambda\tau_x}))
+\esp_\infty(e^{-\lambda\tau^{\{x\}}})\esp_x(e^{-\lambda \tau_x})g_\lambda(x)\,.
\end{equation}
(Remember that  $X(\tau_x)=\infty$ a.s.)
Using Lemma~\ref{lm:enta}, we get that
\begin{equation}
  \label{eq:green}
g_\lambda(x)= \frac { \frac{\lambda\g(x)}{1+\lambda\g(x)} }
                    { \sum_y\frac{\lambda\g(y)}{1+\lambda\g(y)}  }
\end{equation}
We also have the following more general formula.
Let $g_\lambda(x,y)= \lambda\int_0^\infty e^{-\lambda s}\, \P_y(X(s)=x)\, ds$. Then
\beqn
g_\lambda(x,y)=  \frac 1{1+\lambda\g(y)}\, g_\lambda(x)\,.
\eeqn

The last formula describes some correlation function whose definition is motivated by
so-called aging.
\begin{lm}
\label{lm:corage}
Let
\begin{equation}
\label{eq:c}
c_\lambda(\mu)=\int_0^\infty \lambda e^{-\lambda s} ds \int_0^\infty \mu e^{-\mu t} dt\,
\P_\infty(X(u)=X(s)\, \forall u\in[s,s+t])\,.
\end{equation}
Then
\begin{equation}
\label{eq:c1}
c_\lambda(\mu)=\frac {   \sum_x\frac{\lambda\g(x)}{1+\lambda\g(x)}\frac{\mu\g(x)}{1+\mu\g(x)}  }
                     {    \sum_x\frac{\lambda\g(x)}{1+\lambda\g(x)}   }\,.
\end{equation}
\end{lm}
\noindent {\bf Proof}. As for the Green function, we use the Markov property to write that
\beqnn
c_\lambda(\mu)&=& \sum_x  \int_0^\infty \lambda e^{-\lambda s} ds \int_0^\infty \mu e^{-\mu t} dt\,
\P_\infty(X(s)=x\,;\, X(u)=x \, \forall u\in[s,s+t]) \\
&=& \sum_x  \int_0^\infty \lambda e^{-\lambda s} ds \int_0^\infty \mu e^{-\mu t} dt \,
 \P_\infty(X(s)=x)\P_x(\tau_x>t)\\
&=&   \sum_x  g_\lambda(x) (1-\esp_x(e^{-\mu\tau_x})).\quad\qed
\eeqnn

\subsection{Some extension}
\label{ssec:ext}

Let $c>0$ and define the new measure $\g^c=\g+c\delta_\infty$.
The bilinear form $(\E,\D)$ turns out to be also a Dirichlet form when acting on
$L^2(\bN,\g^c)$. The corresponding Markov process can be described as follows:
let $L(t)$ be the local time of $X$ at $\infty$.
($L(t)$ is the unique additive functional whose Revuz measure is $\delta_\infty$.)
Define
\begin{equation}\label{eq:tc}
A^c(t)=t+cL(t) \mbox{ and } X^c(t)=X(A^{-1}(t)).
\end{equation}
Then, under $\P_x$, $X^c$ is a Markov process and its Dirichlet form is
$(\E,\D)$ acting on $L^2(\bN,\g^c)$. Call $\P_x^c$ its law when starting from $x$.

One can then reproduce the same computation as before. In particular we get the expression of
the Green function:
\beqn
g^c_\lambda(x)= \frac { \frac{\lambda\g(x)}{1+\lambda\g(x)} }
                    { c\lambda+\sum_y\frac{\lambda\g(y)}{1+\lambda\g(y)} },\quad x\in\N
                    \eeqn
and, since
$g^c_\lambda(\infty)=1-\sum_x g^c_\lambda(x)$,
\beqn
g^c_\lambda(\infty)= \frac { c\lambda }
                    { c\lambda+\sum_x\frac{\lambda\g(x)}{1+\lambda\g(x)} }.
                    \eeqn
We finally have that
\beqn
g^c_\lambda(x,y)=  \frac 1{1+\lambda\g(y)} g^c_\lambda(x),\quad x,y\in\bN .
\eeqn

\begin{rmk}\label{rmk:un}
It follows from the character of the time change~(\ref{eq:tc})
that for all $c\geq0$, at the entrance time of finite subsets $A$
of $\N$ by $X^c$, starting outside $A$, its distribution is
uniform in $A$. See~\cite{kn:FOT}, Section 6.2, in particular
Theorem 6.2.1.
\end{rmk}



\section{Probabilistic point of view}

\label{sec:prob}

In this section we make an explicit construction for the processes introduced in the
previous section, and study
some of its properties which are relevant for what follows.

Let ${\cal N}=\{(N^{(x)}_t)_{t\geq0},\,x\in\N\}$ be i.i.d.~Poisson
processes of rate 1, with $\s^{(x)}_j$ the $j$-th event time of
$N^{(x)}$, and ${\cal
T}=\{T_0;\,T^{(x)}_i\,i\geq1,x\in\N\}$ be
i.i.d.~exponential random variables of rate 1. ${\cal N}$ and
${\cal T}$ are assumed independent.

For $c\geq0$ and $y\in\bN$, let
\begin{equation}
  \label{eq:tau}
  \Gamma(t)=\Gamma^{c,y}(t)=\g(y)\,T_0+\sum_{x=1}^\infty\g(x)\sum_{i=1}^{N^{(x)}_t}T^{(x)}_i+ct,
\end{equation}
where, by convention, $\sum_{i=1}^{0}T^{(x)}_i=0$ for every $x$.

Let $c\geq0$ be fixed. We define the process $\txcy$ on $\bN$
starting at $y\in\bN$ as follows. For $t\geq0$
\begin{equation}
  \label{eq:txc}
  \txc(t)=\tilde X^{c,y}(t)=\begin{cases}\mbox{}\,y,&\mbox{ if } 0\leq t<\g(y)\,T_0,\\
          \mbox{}\,x,&\mbox{ if } \Gamma(\s^{(x)}_j-)\leq
          t<\Gamma(\s^{(x)}_j)\mbox{ for some }1\leq j<\infty,\\
          \infty,&\mbox{ otherwise}.\end{cases}
\end{equation}

\begin{df}
\label{df:kp} We call $\txcy$ the {\em K-process} with parameters $\g$ and
          $c$. We will also call it sometimes the K($\g,c$)-process for shortness.
\end{df}

\begin{rmk}
\label{rmk:kp}
          The case $c=1$ was introduced by Kolmogorov~\cite{kn:K} as an
          example of a Markov process in a countable state space with an
          instantaneous state. It is known in this context as the {\em first example of
          Kolmogorov} or {\em K1} (Kolmogorov also introduced a second such
          example, known as {\em K2}, which is {\em not} a K-process by our
          definition for any $c\geq0$ and $\g$). The case $c=1$ was then studied in~\cite{kn:KR}
          and~\cite{kn:Ch} (Example 3 in Part II, Chapter 20 of the latter reference),
          where an equivalent construction to the above one is given, and
          elsewhere (e.g.,~\cite{kn:Fr}).
          The general case of $c>0$ isn't really different from the one
          introduced by Kolmogorov; one can go from one case to the other by a
          uniform deterministic time rescaling. The $c=0$ case is
          already considerably different. For one thing, it is not strongly continuous
          (where by strong continuity of a process $Y$ in $\bN$ we mean that
          $\lim_{t\to0}\P_x(Y_t=y)=\delta_{xy}$, the Kronecker's delta, for all $x,y\in\bN$;
          as result of the result of Lemma~\ref{lm:zplae} below, this property is seen to fail
          for the K-process with $c=0$ for $x=y=\infty$), which the K1 process is; following
          L\' evy's classification~\cite{kn:L}, the K-process is of the fourth
          kind for $c=0$, and of the fifth kind for $c>0$.\footnote{According
          to this classification, $\infty$ is termed a {\em fictitious} state,
          rather than an instantaneous state, when $c=0$.} Even though the $c=0$ case is a natural
          extension of the $c>0$ one, we didn't find any explicit mention to it in the literature.
          (In~\cite{kn:L}, though, it is argued in general terms that by looking at a fifth kind
          process outside the instantaneous state, one gets a fourth kind process.)
          Nevertheless we will show
          that precisely the $c=0$ case arises as the scaling limit of a (mean-field)
          disordered spin dynamics (the REM-like trap model in the complete
          graph) at low temperatures. Its irregular behavior near $\infty$,
          associated in particular to its lack of strong continuity, is behind
          the {\em aging phenomenon} exhibited by such dynamics at such
          temperatures~\cite{kn:BD}.
\end{rmk}

\begin{rmk}
\label{rmk:ic} It is clear that $\tilde X^{c,y}(0)=y$ almost
surely for all $y\in\N$. That this also holds for $y=\infty$
follows readily from~(\ref{eq:txc}).
\end{rmk}

\begin{rmk}\label{rmk:cont}
Note on the one hand that $T_0,\,\Gamma(\s^{(x)}_j-),\,\Gamma(\s^{(x)}_j)$
are continuous random variables for every $x\in\bN$ and $j\geq1$, and
on the other hand that $\tilde X^{c,y}$ is almost surely continuous off
$\{\g(y)\,T_0;\Gamma(\s^{(x)}_j-),\Gamma(\s^{(x)}_j),\,x\in\bN,j\geq1\}$.
These readily imply that every $s\geq0$ is almost surely a continuity
point of $\txcy$.
\end{rmk}

\begin{rmk}\label{rmk:yinf}
It readily follows from~(\ref{eq:txc}) that
\begin{equation}\label{eq:yi}
    \txcy(t)= \txci(t-\g_y\,T_0),\quad\mbox{for } t\geq\g_y\,T_0.
\end{equation}
\end{rmk}

\begin{prop}
\label{prop:mark} $\txc$ is cadlag and Markovian.
\end{prop}

\begin{rmk}\label{rmk:krc}
A treatment of the case $c=1$ can be found in~\cite{kn:KR} and~\cite{kn:Ch}. Even though
both have a construction equivalent to ours, complete proofs of some key properties
of the constructed process, like the Markov one, are not presented. For this
reason, and in order to include the $c=0$ case as well, we present below a
proof of Proposition~\ref{prop:mark}.
\end{rmk}

The proof is based on strongly approximating $\txc$ in Skorohod
space by Markov processes that we now define. For
$n\geq1$ and $y\in\{1,\ldots,n,\infty\}$, let
\begin{equation}
  \label{eq:taun}
  \Gamma_n(t)=\Gamma_n^{c,y}(t)=\g(y)\,T_0+\sum_{x=1}^n\g(x)\sum_{i=1}^{N^{(x)}_t}T^{(x)}_i+ct
\end{equation}
and
\begin{equation}
  \label{eq:txcn}
  \tilde X_n^{c,y}(t)=\begin{cases}\mbox{}\,y,&\mbox{ if } 0\leq t<\g(y)\,T_0,\\
                         \mbox{}\,x,&\mbox{ if } \Gamma_n(\s^{(x)}_j-)\leq
  t<\Gamma_n(\s^{(x)}_j)\mbox{ for some }1\leq x\leq n,\,j\geq1,\\\infty,&\mbox{ otherwise}.\end{cases}
\end{equation}
\begin{rmk}
\label{rmk:xn1} We note that $\tilde X_n^{0,y}$ never visits
$\infty$, even when $y=\infty$. See next remark.
\end{rmk}
\begin{rmk}
\label{rmk:xn2} The order in which the sites of $\{1,\ldots,n\}$
are visited by $\tilde X_n^{c,y}$ (in case $y$ is finite, after
leaving the initial state) is given by the respective
(chronological) order of $\{\s^{(x)}_j;\,1\leq x\leq
n,\,j\geq1\}$. Let us denote the latter set by
$\S^n=\{S^n_1,S^n_2,\ldots\}$, with $S^n_1<S^n_2<\ldots$. Then
$\S^n$ is a Poisson point process of rate $n$, each point of which
is labeled according to a different element of an i.i.d.~family of
uniform in $\{1,\ldots,n\}$ random variables. This implies that
the jump probabilities of $\tilde X_n^{c,y}$ from any site in case
$c=0$, and from $\infty$ in case $c>0$, are uniform in
$\{1,\ldots,n\}$, and also implies that $\tilde X_n^{0,\infty}(0)$
is uniformly distributed in $\{1,\ldots,n\}$ (since it is the
label of $S^n_1$; see previous remark).

In case $c>0$, $c(S^n_i-S^n_{i-1})$, $i\geq1$, where
$S^n_0\equiv0$, represent the successive holding times at
$\infty$. It is clear then that these times form an
i.i.d.~sequence of exponential random variables of mean $c/n$.
\end{rmk}
We have the following two results.
\begin{lm}
\label{lm:mark} $\txc_n$ is cadlag and Markovian for every
$n\geq1$ and $y\in\{1,\ldots,n,\infty\}$.
\end{lm}
\begin{lm}
\label{lm:approx}  $\txc_n\to\txc$ as $n\to\infty$ almost surely
in the Skorohod norm for every $y\in\bN$.
\end{lm}

\noindent{\bf Proof of the first assertion of
Proposition~\ref{prop:mark}}

The first assertion of Lemma~\ref{lm:mark} and
Lemma~\ref{lm:approx} readily establish the first assertion of
Proposition~\ref{prop:mark} (see~\cite{kn:EK}). \qed

\bigskip

\noindent{\bf Proof of Lemma~\ref{lm:mark}}

Let the starting point $y$ be fixed.

For $c=0$, $\txc_n$ is the following Markov process on
$\{1,\ldots,n\}$. $\tilde X_n^{c,y}$ starts at $y$ if
$y\in\{1,\ldots,n\}$; $\tilde X_n^{c,\infty}$ has uniform initial
distribution. When at $x\in\{1,\ldots,n\}$, it waits an
exponential time of mean $\g(x)$ and then jumps uniformly at
random to a site in $\{1,\ldots,n\}$ (which could be $x$ again).
See Remarks~\ref{rmk:xn1} and~\ref{rmk:xn2} above.

For $c>0$, $\txc_n$ is the following Markov process on
$\{1,\ldots,n,\infty\}$. $\tilde X_n^{c,y}$ starts at $y$. When at
$x\in\{1,\ldots,n\}$, it waits an exponential time of mean $\g(x)$
and then jumps deterministically to $\infty$. When at $\infty$, it
waits an exponential time of mean $c/n$, and then jumps uniformly
at random to a site in $\{1,\ldots,n\}$. See Remark~\ref{rmk:xn2}
above. \qed

\bigskip

\noindent{\bf Proof of Lemma~\ref{lm:approx}}

Let $y$ be fixed, and suppose $n\geq y$ if $y\in\N$.
We show the almost sure validity of $(c)$ of
Proposition 5.3 in Chapter 3 of~\cite{kn:EK} (page~119).

For $m\in\N$, let $\d_m=\mbox{diam}\{x\in\bN:\,x>m\}=(m+1)^{-1}$ and
$\{S_1^m<S_2^m<\ldots\}=\{\s^{(x)}_j,\,j\geq1,1\leq x\leq m\}$,
with the latter being well defined almost surely. Fix $T>0$
and let
\begin{equation}\label{lmn}
L^m_n=\min\{i\geq1:\,\Gamma_n(S^m_i)\geq T\},
\end{equation}
which is almost surely finite, and make $S^m_0\equiv0$.
Notice that $L^m_n=L^m_n(y)$, is nondecreasing in $y$,
where $y$ is the starting point, and thus
\begin{equation}\label{lmn1}
\max_{y\in\bN}L^m_n(y)=L^m_n(\infty).
\end{equation}
We can now almost surely find $n_m$ so large that $\min_{0\leq i\leq
L^m_n-1}[\Gamma_n(S^m_{i+1}-)-\Gamma_n(S^m_{i})]>0$ for $n\geq
n_m$\footnote{We can take $n_m\equiv1$ when $c>0$.}. For these $n$
then define $\l_n^m:[0,\Gamma_n(S^m_{L^m_n})]\to\R^+$ inductively as
follows.
\begin{equation}\label{ln0}
    \l_n^m(t)=t,\mbox{ if } 0\leq t<\g(y)\,T_0,
\end{equation}
and for $0\leq i\leq L^m_n-1$ and $\Gamma_n(S^m_i)\leq
t\leq\Gamma_n(S^m_{i+1})$, let
\begin{equation}
  \label{eq:ln}
  \l_n^m(t)=\begin{cases}\Gamma(S^m_i)+
  \frac{\Gamma(S^m_{i+1}-)-\Gamma(S^m_{i})}{\Gamma_n(S^m_{i+1}-)-\Gamma_n(S^m_{i})}\,[t-\Gamma_n(S^m_{i})],&
  \mbox{ if } \Gamma_n(S^m_{i})\leq t\leq\Gamma_n(S^m_{i+1}-),\\\mbox{}\quad\quad
  \Gamma(S^m_{i+1}-)-\Gamma_n(S^m_{i+1}-)+t,&\mbox{ if } \Gamma_n(S^m_{i+1}-)\leq t\leq\Gamma_n(S^m_{i+1}).
\end{cases}
\end{equation}
It has the following properties. For all $T>0$, $m\in\N$ and $n\geq m\vee n_m\vee y$
\begin{eqnarray}
  \label{eq:ln1/2}\l_n^m(t)&\geq&t,\,\,0\leq t\leq T,\\
  \label{eq:ln1}
  \sup_{0\leq t\leq T}|\l_n^m(t)-t|
  &\leq&\max\left\{
  \Gamma(S^m_{i+1}-)-\Gamma_n(S^m_{i+1}-);\,0\leq i\leq
  L^m_n(\infty)-1\right\}
\end{eqnarray}
(where we have made use of~(\ref{lmn1})),
and
\begin{equation}
  \label{cas}
\mbox{the right hand side of~(\ref{eq:ln1}) vanishes almost surely
as $n\to\infty$.}
\end{equation}
Furthermore,
\begin{equation}
  \label{eq:ln2}
  \sup_{0\leq t\leq
  T}\mbox{dist}\left(\txc(\l_n^m(t)),\txc_n(t)\right)\leq\d_m,
\end{equation}
since for $t\in[0,T]$, $\txc(\l_n^m(t))$ and $\txc_n(t)$ coincide
when either one is in $\{1,\ldots,m\}$.

{}From~(\ref{eq:ln1}), for every $m\in\N$ there almost surely exists $n'_m\geq
n_m$ such that for $n\geq n'_m$
\begin{equation}
  \label{eq:ln3}
  \sup_{0\leq t\leq T}|\l_n^m(t)-t|\leq\d_m
\end{equation}
and~(\ref{eq:ln2}) hold. We may assume $(n'_m)$ is strictly
increasing.

For $n\geq n'_1$, let $m_n=i$ when $n'_i\leq n<n'_{i+1}$ and
$\tilde \l_n=\l_n^{m_n}$. We then have for $T>0$
\begin{eqnarray}
  \label{eq:ln4}
  &\sup_{0\leq t\leq
  T}\mbox{dist}\left(\txc(\tilde\l_n(t)),\txc_n(t)\right)\to0&\\
  \label{eq:ln5}
  &\sup_{0\leq t\leq T}|\tilde\l_n(t)-t|\to0&
\end{eqnarray}
almost surely as $n\to\infty$, and the above mentioned condition $(c)$ is
verified. \qed

\bigskip

We will also need the further lemma to prove the second assertion
of Proposition~\ref{prop:mark}.

\begin{lm}
\label{lm:ytoinf} For every $t\geq0$, $\txcy(t)\to\txci(t)$ as $y\to\infty$ almost surely.
\end{lm}

\noindent{\bf Proof}

The case $t=0$ is clear. For $t>0$, since we are taking $y\to\infty$,
we may assume that $\g(y)\,T_0\leq t$ and then from~(\ref{eq:yi}) we
have that
\begin{equation}\label{eq:yi1}
    |\txcy(t)-\txci(t)|=|\txci(t-\g(y)\,T_0)-\txci(t)|
\end{equation}
and the result follows from Remark~\ref{rmk:cont}.
 \qed

\bigskip

\noindent{\bf Proof of the second assertion of
Proposition~\ref{prop:mark}}

Lemma~\ref{lm:mark} implies that for arbitrary $m\geq1$, $0\leq
t_1<\ldots<t_{m+1}$ and bounded continuous functions $f_1,\ldots,
f_{m+1}$, we have that
\begin{eqnarray}\nn
    &\esp\left[f_1(\txc_n(t_1))\ldots
    f_m(\txc_n(t_m))\,f_{m+1}(\txc_n(t_{m+1}))\right]&\\
    \label{eq:mk1}
    &=\esp\left[f_1(\txc_n(t_1))\ldots
    f_m(\txc_n(t_m))\,\Psi^n_{t_{m+1}-t_{m}}f_{m+1}(\txc_n(t_{m}))\right],&
\end{eqnarray}
where $\Psi^n$ is the semigroup of $\txc_n$, i.e., for $t\geq0$, a continuous
function $f$ and $y\in\bN$
\begin{equation}\label{eq:psin}
    \Psi^n_tf(y)=\esp\left[f(\txcy_n(t))\right].
\end{equation}
By Lemma~\ref{lm:approx}, the left hand side of~(\ref{eq:mk1})
converges to
\begin{equation}\label{eq:mk2}
    \esp\left[f_1(\txc(t_1))\ldots
    f_m(\txc(t_m))\,f_{m+1}(\txc(t_{m+1}))\right]
\end{equation}
as $n\to\infty$. Let us estimated the right hand side
of~(\ref{eq:mk1}) by
\begin{equation}\label{eq:mk3}
    \esp\left[f_1(\txc_n(t_1))\ldots
    f_m(\txc_n(t_m))\,\Psi_{t_{m+1}-t_{m}}f_{m+1}(\txc_n(t_{m}))\right]+\e_n,
\end{equation}
where for $t\geq0$, a continuous function $f$ and $y\in\bN$
\begin{equation}\label{eq:psi}
    \Psi_tf(y)=\esp\left[f(\txcy(t))\right],
\end{equation}
and
\begin{equation}\label{eq:en}
|\e_n|\leq\mbox{const }\sup_y|\Psi^n_{t_{m+1}-t_{m}}f_{m+1}(y)-\Psi_{t_{m+1}-t_{m}}f_{m+1}(y)|.
\end{equation}

{From} Lemma~\ref{lm:ytoinf}, we have that $\Psi_{t_{m+1}-t_{m}}f_{m+1}(\cdot)$ is continuous, and
now Lemma~\ref{lm:approx} implies that the left term
of~(\ref{eq:mk3}) converges to
\begin{equation}\label{eq:mk4}
    \esp\left[f_1(\txc(t_1))\ldots
    f_m(\txc(t_m))\,\Psi_{t_{m+1}-t_{m}}f_{m+1}(\txc(t_{m}))\right]
\end{equation}
as $n\to\infty$.

Let us now examine the right hand side of~(\ref{eq:en}). We
first relabel $f_{m+1}=g$, $t_{m+1}-t_m=s$ and $\g(y)=\g_y$. We
have that
\begin{eqnarray}\nn
&\Psi^n_{t_{m+1}-t_{m}}f_{m+1}(y)-\Psi_{t_{m+1}-t_{m}}f_{m+1}(y)
  =\esp\left[g(\txcy_n(s))-g(\txcy(s))\right]&\\\label{eq:en1}
&=\esp\left[g(\txcy_n(s))-g(\txcy(\tl_n(s)))\right]+\esp\left[g(\txcy(\tl_n(s)))-g(\txcy(s))\right],&
\end{eqnarray}
with $\tl_n$ as defined in the paragraph of~(\ref{eq:ln4}), with
$T>s$. From~(\ref{eq:ln2}), it follows that the sup in $y$ of the
absolute value of the first expected value in the right hand side
of~(\ref{eq:en1}) vanishes as $n\to\infty$ (since $g$ is
continuous, and thus uniformly continuous since $\bN$ is compact).

Lemma~\ref{lm:approx} now implies that there exists a sequence $k_n$ going to
infinity as $n\to\infty$ such that as $n\to\infty$
\begin{equation}\label{eq:en2}
\max_{1\leq y\leq k_n}
\left|\esp\left[g(\txcy(\tl_n(s)))-g(\txcy(s-\g_y\,T_0))\right]\right|
\to 0.
\end{equation}

We now note that, from~(\ref{eq:ln1/2},\ref{cas}), $\tln(s)\geq s$, and that $\tln(s)\to s$
as $n\to\infty$ uniformly in $y$ almost surely. From this and~(\ref{eq:yi}) we then have
\begin{eqnarray}\nn
&&\sup_{y>k_n}
  \left|\esp\left[g(\txcy(\tl_n(s)))-g(\txcy(s-\g_y\,T_0))\right]\right|\\
&\leq&\nn
\esp\left[\sup_{y>k_n}\left|g(\txci(\tl_n(s)-\g_y\,T_0))-g(\txci(s-\g_y\,T_0))\right|;\,
       \g_y\,T_0<s\right]\\\label{eq:en4}
\+2\|g\|\sup_{y>k_n}\P(\g_y\,T_0\geq s).
\end{eqnarray}
It is clear that the latter summand in the right hand side of~(\ref{eq:en4}) vanishes
as $n\to\infty$. And so does the former one, since
$s$ is almost surely a continuity point of $\txci$ (see Remark~\ref{rmk:cont} above).

We have thus concluded that $|\e_n|\to0$ as $n\to\infty$,
and then from~(\ref{eq:mk1}-\ref{eq:mk4}) we have that
\begin{eqnarray}\nn
    &\esp\left[f_1(\txc(t_1))\ldots
    f_m(\txc(t_m))\,f_{m+1}(\txc(t_{m+1}))\right]&\\\label{eq:mki}&=
    \esp\left[f_1(\txc(t_1))\ldots
    f_m(\txc(t_m))\,\Psi_{t_{m+1}-t_{m}}f_{m+1}(\txc(t_{m}))\right],&
\end{eqnarray}
and the Markov property is established.
\qed

\begin{rmk}
\label{rmk:sg} $\Psi$ defined in~(\ref{eq:psi})
is the semigroup of $\txc$.
\end{rmk}

\medskip

Next follows a result establishing in particular the lack of strong
continuity of the K-process with $c=0$.

\begin{lm}
\label{lm:zplae} For every $y\in\bN$, we have that $\P(\tilde
X^{0,y}(t)=\infty)=0$ for  every $t>0$.
\end{lm}

\begin{rmk}
\label{rmk:posp} The statement of Lemma~\ref{lm:zplae} does not
hold for $c>0$. In this case, it can actually be shown that
$\P(\txc(t)=\infty)>0$ for every $t>0,\,y\in\bN$. It can also be
shown that the process is strongly continuous in this case.
\end{rmk}

\begin{rmk}
\label{rmk:no} We note that $\tilde X^{c,\infty}(0)=\infty$ almost
surely for every $c\geq0$.
\end{rmk}

\noindent{\bf Proof of Lemma~\ref{lm:zplae}}

For $m\geq1$, and $t>0$, let $\tmt$ be the time spent by $\tilde
X^{0,y}$ outside $\{1,\ldots,m\}$ up to time $t$.
Clearly
\begin{equation}\label{eq:tmt}
    \tmt=\sum_{x>m}\g(x)\sum_{i=1}^{N^{(x)}_{\Xi_t}}T^{(x)}_i,
\end{equation}
where $\Xi$ is the inverse function of $\Gamma$. It is also  clear
that
\begin{equation}\label{eq:zp1}
    \int_0^t1_{\{\infty\}}(\tilde X^{0,y}(s))\,ds\leq\tmt
\end{equation}
for every $m\geq1$ and $t>0$, where $1_{\cdot}$ is the usual
indicator function, and that
\begin{equation}\label{eq:zp2}
   \tmt\to0
\end{equation}
almost surely as $m\to\infty$ for every $t>0$. Thus the left hand
side of~(\ref{eq:zp1}) vanishes almost surely and dominated
convergence implies that
\begin{equation}\label{eq:zp3}
    \esp\left(\int_0^t1_{\{\infty\}}(\tilde X^{0,y}(s))\,ds\right)=
    \int_0^t\P(\tilde X^{0,y}(s)=\infty)\,ds=0
\end{equation}
for every $t$. This proves the assertion of the lemma for
Lebesgue-almost every $t$. The Markov property of $\tilde X^{0,y}$
can now be used to extend the result to every $t$.
 \qed

\bigskip

We close this section with a computation related to the Green
function of $\tilde X^c$; this will lead to an identification of
$\tilde X^c$ above and $X^c$ defined in Subsection~\ref{ssec:ext}.

Let $\tau^{\{x\}}=\inf\{t\, ;\, \tilde X^c(t)=x\}$. We have that
under $\P_\infty$
\begin{equation}
  \label{eq:ag70}
  \tau^{\{x\}}=\gx_c(\sigma^{(x)}_1),
\end{equation}
where for $x\in\N,\,s\geq0$
\begin{equation}
  \label{eq:ag7}
  \gx_c(s)=\sum_{y\ne x}\g_y\sum_{i=1}^{N^{(y)}_s}T^{(y)}_i+cs.
\end{equation}
It is now straightforward to compute the Laplace transform of
$\tau^{\{x\}}$ for the process started at $\infty$.  We obtain
\begin{eqnarray}\nonumber
  &\esp_\infty\!\!\left(e^{-\lambda\tau^{\{x\}}}\right)=
  \esp_\infty\!\!\left(e^{-\lambda\gx_c(\sigma^{(x)}_1)}\right)=&\\
    \label{eq:ag71}
  &\int_0^\infty\esp\!\left(e^{-\lambda\gx_c(s)}\right)e^{-s}\,ds=\int_0^\infty\esp\!
  \left(\exp\left\{-\lambda\sum_{y\ne x}\g_y
  \sum_{i=1}^{N^{(y)}_s}T^{(y)}_i\right\}\right)
  e^{-(1+c)s}\,ds,&
\end{eqnarray}
for $\lambda>0$, where in the second equality we have used the independence of
$\sigma^{(x)}_1$ and (the random variables in)
$\sum_{y\ne x}\g_y\sum_{i=1}^{N^{(y)}_s}T^{(y)}_i$.
We leave it as an exercise to compute the expectation inside the integral
in~(\ref{eq:ag71}), and to conclude that the integral equals
\begin{equation}
  \label{eq:ag71a}
\left(1+c\lambda+\sum_{y\ne x}\frac{\lambda\g_y}{1+\lambda\g_y}\right)^{-1}.
\end{equation}
We note that this expression is the same as that for the
corresponding transform for $X^c$ in Section~\ref{sec:anal}.

Now, since the only transitions are from states in $\N$ to
$\infty$ and back, we have a decomposition as in~(\ref{eq:marpro})
for the Green kernel of $\tilde X^c$ starting at $\infty$, 
and we get the case of general initial condition from the case of $\infty$
initial condition as in the computation in Section~\ref{sec:anal}. 
We readily conclude from the remark at the end 
of the previous paragraph that the Green
functions of $\tilde X^c$ and $X^c$ coincide, and thus, since
these are both c\`adl\`ag Markov processes, they must have the
same distribution for any initial law.



\section{A characterization result}
\label{sec:char}

The striking property of K-processes that at the entrance time of the process
in finite subsets (starting from outside) the distribution is uniform
(see Remark~\ref{rmk:un}) leads to a natural question: which other processes
have this property? Below we see that, under natural assumptions, the answer is
{\em none}, that is, that property characterizes K-processes.

\begin{theo}
\label{teo:char}
Let $\g$ be as in~(\ref{eq:sum}) and $Y=(Y(t),t\geq 0)$ be a
process on $\bN$ with the following four properties.
\begin{itemize}
\item[(i)] Y is c\`adl\`ag.

\item[(ii)] Y is strong Markov.

\item[(iii)] Starting from any point $i\in\N$, $Y$ waits for an
exponential time of mean $\g(i)$ before jumping.

\item[(iv)] Starting from $\infty$, for any finite $A\subset\N$, we have
$\tau_{A}<\infty$ almost surely, where
\begin{equation}
 \tau_{A}=\inf\{t\geq0:\,Y(t)\in A\},
\end{equation}
with $\inf\emptyset=\infty$, and the law of $Y(\tau_{A})$ is uniform on $A$.
\end{itemize}
Then, $Y$ is a K-process with parameters $\g$ and $c$, for some $c\geq0$.
\end{theo}

\begin{rmk}\label{fuku}
Fukushima and collaborators, as well as other authors, have
recently studied one-point extensions of certain Markov processes
beyond a killing time (see~\cite{kn:FT,kn:CFY} and references
therein). The K-processes can be viewed as one-point extensions of
processes in $\N$ that are killed after the first jump. With this
point of view, and even though the K-processes don't satisfy some
of the conditions in the above references (like Condition (A.2)
in~\cite{kn:CFY}; another condition would require $c=0$ in our
case), Theorem~\ref{teo:char} is similar (in its particular
context) to their results. But there is an important difference in
that, while they depart from a reversibility condition (more
generally, a duality condition) with respect to an excessive
measure for the process, we have a condition on the jump rates and
entrance laws. It is nevertheless remarkable that entrance laws
play a crucial role in their approach (it also could be said that
for the K-processes the jump rates are directly related to a
stationary measure for the process).
\end{rmk}

\noindent{\bf Proof of Theorem~\ref{teo:char}} The strategy is to
consider the process restricted to $\{1,\ldots,n,\infty\}$, and
show that it must have the same distribution as $\txc_n$. The
result then follows by taking $n\to\infty$.

We start by showing that from any state in $\N$, $Y$ jumps to
$\infty$ almost surely. Let $i,j\in\N$ be such that $i\ne j$, and
for $n\geq i\vee j$ let
\begin{equation}\label{eq:an}
A_n=\{1,\ldots,n\}.
\end{equation}
Then
\begin{equation}
\{Y(\tau_{A_{n+1}})=j\}\,\cup\,\{Y(\tau_{A_{n+1}})=i,\,Y(\tau'_{n+1})=j\}\,
 \subset\,\{Y(\tau_{A_{n+1}\setminus\{i\}})=j\},
\end{equation}
where
\begin{equation}
\tau_n'=\inf\{t\geq\tau_{A_{n}}:\,Y(t)\ne Y(\tau_{A_{n+1}})\}.
\end{equation}
Thus,
\begin{equation}
\P_\infty(Y(\tau_{A_{n+1}})=j)+\P_\infty(Y(\tau_{A_{n+1}})=i,\,Y(\tau'_{n+1})=j)\leq
\P_\infty(Y(\tau_{A_{n+1}\setminus\{i\}})=j),
\end{equation}
and using (ii-iv)
\begin{equation}
\frac1{n+1}+\frac1{n+1}\,p_{ij}\leq \frac1{n},
\end{equation}
where $p_{ij}$ is the transition probability from $i$ to $j$. It
follows that $p_{ij}\leq1/n$, and since $n$ can be taken
arbitrarily large, we conclude that $p_{ij}=0$, and the claim at
the beginning of the paragraph follows.

Now let us consider the process obtained from $Y$ by suppressing
jumps outside $\bar A_n:=A_n\cup\{\infty\}$ (see~(\ref{eq:an})).
Let us call it $Y_n$. More precisely let
\beqnn
\A_n(t)=\int_0^t\1_{\bar A_n}(Y(s))\, ds \hbox{ \ and \ } Y_n(t)=Y(\A_n^{-1}(t)),
\eeqnn
where $\1_{A}$ is the usual indicator function of a set $A$, and $\A_n^{-1}$
is the right-continuous inverse of $\A_n$.

It is readily seen that $\A_n(t)\uparrow t$ as $n\uparrow\infty$ uniformly
in $t\leq T$ for every $T$. This and (i) implies that
\begin{equation}\label{eq:yntoy}
Y_n\to Y\,\,\mbox{ as $n\to\infty$ in Skorohod space.}
\end{equation}
$Y_n$ also satisfies (i) and (ii) (see~\cite{kn:S}, Theorem 65.9).
Starting at $i\leq n$, it waits an exponential time of mean
$\g(i)$ and then jumps.

The state space of $Y_n$ may be either $\bar A_n$ or $A_n$. (The
latter possibility happens if $Y$ is the K($\g,0$)-process; note
that $\int_0^\infty\1_{\{\infty\}}(Y(s))\, ds=0$ almost surely in
that case --- see the proof of Lemma~\ref{lm:zplae} for an
argument.) If the latter case happens, then $Y_n$ is a continuous
time Markov chain on $A_n$ satisfying (i-iii). To completely
characterize it, we need only determine the transition
probabilities. But Property (iv) of $Y$ implies that these must be
uniform, that is, $p_{ij}\equiv1/n$. This means that $Y_n$ is
equidistributed with $\tilde X^0_n$ defined in~(\ref{eq:txcn}).
Now this, Lemma~\ref{lm:approx} and~(\ref{eq:yntoy}) imply that
$Y$ is equidistributed with $\tilde X^0$.

I remains to consider the case where the state space of $Y_n$ is $\bar A_n$.
In this case $Y_n$ clearly also satisfies (iii-iv). We need only
determine the mean holding time at $\infty$, say $\g_n(\infty)$. For that we reason
as follows.

We can obtain $Y_n$ by suppressing jumps of $Y_{n+1}$ outside
$\{1,...,n\}\cup\{\infty\}$. Since upon leaving $\infty$, the process $Y_{n+1}$ has
probability $1/(n+1)$ to jump to $n+1$, we see that the holding time at
$\infty$ in $Y_{n}$ can be seen as a sum of independent holding times at
$\infty$ in $Y_{n+1}$. The number of terms in the sum is a geometric random
variable with success parameter $1/(n+1)$, independent of the holding times
at $\infty$ in $Y_{n+1}$. We conclude that
\beqnn
\g_n(\infty)=\g_{n+1}(\infty)\,\frac {n+1}n,
\eeqnn
and thus that
$\g_n(\infty)=c/n$ for some constant $c\geq0$. We then see that $Y_n$ is
equidistributed with $\tilde X^c_n$ defined in~(\ref{eq:txcn}) for every
$n\geq1$, and the conclusion that $Y$ is equidistributed with $\tilde X^c$
follows exactly as above. \qed



\section{A scaling limit for the REM-like trap model}
\label{sec:bou}

The REM-like trap model~\cite{kn:BD} can be described as a continuous time
symmetric Markov chain $Y_n=\{Y_n(t),\,t\geq0\}$ in the complete graph $K_n$
with $n$ vertices such that the average holding times $\tau:=\{\tau_x,\,x\in K_n\}$
is an i.i.d.~family of positive random variables equidistributed with a r.v.~$\tau_0$
which is in the basin of attraction of a stable law of degree $\a<1$, that is
\begin{equation}
  \label{eq:tau0}
  \P(\tau_0>t)=\frac{L(t)}{t^\a},\,t>0,
\end{equation}
where $L$ is a slowly varying function at infinity.

We will show in this section that in an appropriate sense, in an appropriate
time scale, $Y_n$ converges in distribution as $n\to\infty$ to a K-process
with $c=0$.

We start by identifying the vertices of $K_n$ with $A_n=\{1,\ldots,n\}$ for all
$n\geq1$, in such a way that $\{\tau^{(n)}_i,\,1\leq i\leq n\}$ is in decreasing
order (that is, $(\tau^{(n)}_1,\ldots,\tau^{(n)}_n)$ is the reverse order
statistics of $(\tau_1,\ldots,\tau_n)$, an i.i.d.~sample of size $n$ of
$\tau_0$.

We can describe $Y_n$ then as a continuous time Markov chain in $A_n$ with
mean holding time at $i\in A_n$ given by $\tau^{(n)}_i$ and uniform in $A_n$ transition
probabilities for all starting point $i\in A_n$.

Let us view $(\tau^{(n)}_1,\ldots,\tau^{(n)}_n)$ as a random measure $\g_n$ on
$\N$ such that
\begin{equation}
  \label{eq:gn}
  \g_n(\{i\})=\begin{cases}\tau^{(n)}_i,&\,\mbox{ if } i\in A_n,\\
                         0,&\,\mbox{ otherwise}.
          \end{cases}
\end{equation}

We refer the reader to Section 3 of~\cite{kn:FIN} for more on the context
of the next result. We present the main points below.

Consider the increasing L\' evy process
$V_x,\,x\in\R,\,V_0=0$, with stationary
and independent increments given by
\begin{equation}
\label{eq:gf}
\esp\left[e^{ir(V_{x+x_0}-V_{x_0})}\right]=
e^{\a x\int_0^\infty(e^{irw}-1)\,w^{-1-\a}\,dw}
\end{equation}
for any $x_0 \in \R$ and $x \geq 0$.
Let $\rho$ be the (random) Lebesgue-Stieltjes measure on the
Borel sets of $\R$ associated to $V$, i.e.,
\begin{equation}
\label{eq:mu}
\rho((a,b])=V_b-V_a,\,a,b\in\R,\, a<b.
\end{equation}
Then
\begin{equation}
\label{eq:dmu} d\rho=dV=\sum_j\,w_j\,\d(x-x_j),
\end{equation}
where the (countable) sum is over the indices of an inhomogeneous
Poisson point process $\{(x_j,w_j)\}$ on $\R\times(0,\infty)$ with
density $dx\,\a w^{-1-\a}\,dw.$

Let now
\begin{equation}
\label{eq:c_n}
c_n=\left(\inf\{t\geq0:\P(\tau_0>t)\leq n^{-1}\}\right)^{-1}.
\end{equation}
and $\tilde\g_n=c_n\g_n$, that is,
$\tilde\g_n$ is a (random) measure in $\N$ such that
\begin{equation}
  \label{eq:tgn}
  \tilde\g_n(\{i\})=\begin{cases}c_n\tau^{(n)}_i,&\,\mbox{ if } i\in A_n,\\
                         0,&\,\mbox{ otherwise}.
          \end{cases}
\end{equation}
Let $\g=\{\g(i),\,i\in\N\}$ denote the weights of $\rho$
in $[0,1]$ in decreasing order, that is, making $\RR=\{\g(\{x\}),\,x\in[0,1]\}$,
  \begin{equation}
    \label{eq:ga}
    \g(1)=\max\RR;\,\,\,\g(i)=\max\left[\RR\setminus\{\g(1),\ldots,\g(i-1)\}\right],\,i\geq2.
  \end{equation}

\begin{rmk}\label{rmk:ga}
$\g$ thus defined almost surely satisfies the conditions on the paragraph of~(\ref{eq:sum}).
\end{rmk}

\begin{theo}
  \label{teo:scabou}
  Let $\tilde Y_n$ be the process in $A_n$ such that for $t\geq0$, $\tilde
  Y_n(t)=Y_n(c_n^{-1}t)$. Suppose $\tilde Y_n(0)\equiv{\cal Y}_n$ converges
  weakly to a random variable $\cal Y$ in $\bN$. Then, as $n\to\infty$,
  \begin{equation}
    \label{eq:sb}
    (\tilde Y_n,\tilde\g_n)\Rightarrow(Y,\g),
  \end{equation}
where, given $\g$, $Y$ is a K($\g,0$)-process with $Y(0)$
distributed as $\cal Y$, and $\Rightarrow$ denotes weak
convergence in the product of the Skorohod topology and the vague
topology in the space of finite measures on $\bN$.
\end{theo}

\noindent{\bf Proof} We may assume that ${\cal Y}_n\to\cal Y$ as
$n\to\infty$ almost surely. Following the strategy in Section 3
of~\cite{kn:FIN}, we will couple $(\tilde Y_n,\tilde\g_n)$ to
$(Y,\g)$ and establish~(\ref{eq:sb}) as a strong convergence.

For $i\in A_n$, let
\begin{equation}
\label{eq:tn}
\bar\tau^{(n)}_i=\frac1c_n\,g_n\!\left(V_{i/n}-V_{(i-1)/n}\right),
\end{equation}
where $g_n$ is defined as follows. Let $G:[0,\infty)\to[0,\infty)$
satisfy
\beq
\label{eq:coup6}
\P(V_1>G(x))=\P(\tau_0>x)\quad\mbox{for all $x\geq 0$}
\eeq
and let $g_n:[0,\infty)\to[0,\infty)$
be defined as
\beq
\label{eq:coup7}
g_n(x)=c_n G^{-1}\!\left(n^{1/\a}x\right)\quad\mbox{for all $x\geq 0$.}
\eeq

Let $\{\hat\tau^{(n)}_i,\,i\in A_n\}$ be
$\{\bar\tau^{(n)}_i,\,i\in A_n\}$ in decreasing order, and
\begin{equation}
  \label{eq:hgn}
  \hat\g_n(\{i\})=\begin{cases}c_n\hat\tau^{(n)}_{i},&\,\mbox{ if } i\in A_n\\
                         0,&\,\mbox{ otherwise}.
          \end{cases}
\end{equation}
It readily follows from Proposition 3.1 in~\cite{kn:FIN} that
$\hat\g_n$ and $\tilde\g_n$ have the same distribution for every $n$,
and that almost surely
\begin{equation}
  \label{eq:sb1}
\hat\g_n\to\g\,\mbox{ as } n\to\infty
\end{equation}
(where the first $\to$ in~(\ref{eq:sb1}) means vague convergence).

Let now $\n$ and $\T$ be as in Section~\ref{sec:prob}. For $t\geq0$, let
\begin{eqnarray}
  \label{eq:htaun}
  \hat\Gamma_n(t)\=\hat\g_n({\cal Y}_n)\,T_0+\sum_{x=1}^n\hat\g_n(x)
  \sum_{i=1}^{N^{(x)}_t}T^{(x)}_i,\\
  \label{eq:htau}
  \hat\Gamma(t)\=\g({\cal Y})\,T_0+\sum_{x=1}^\infty\g(x)\sum_{i=1}^{N^{(x)}_t}T^{(x)}_i,
\end{eqnarray}
where we write $\hat\g_n(x)$ and $\g(x)$ for $\hat\g_n(\{x\})$ and
$\g(\{x\})$, respectively, and
\begin{eqnarray}
  \label{eq:hyn}
  \hat Y_n(t)\=\begin{cases}\mbox{}\,{\cal Y}_n,&\mbox{ if } 0\leq t<\hat\g_n({\cal Y}_n)\,T_0,\\
                         \mbox{}\,x,&\mbox{ if } \hat\Gamma_n(\s^{(x)}_j-)\leq
  t<\hat\Gamma_n(\s^{(x)}_j)\mbox{ for some }1\leq x\leq n,\,j\geq1.\end{cases}\\
  \label{eq:hy}
  \hat Y(t)\=\begin{cases}\mbox{}\,{\cal Y},&\mbox{ if } 0\leq t<\g({\cal Y})\,T_0,\\
                         \mbox{}\,x,&\mbox{ if } \hat\Gamma(\s^{(x)}_j-)\leq
  t<\hat\Gamma(\s^{(x)}_j)\mbox{ for some }1\leq x<\infty,\,j\geq1\\\infty,&\mbox{
                         otherwise}.\end{cases}
\end{eqnarray}
See~(\ref{eq:tau}),\,\ref{eq:txc}),~(\ref{eq:taun},\,\ref{eq:txcn}) above. One
readily checks that $(\hat Y_n,\hat\g_n)$ has the same distribution as
$(\tilde Y_n,\tilde\g_n)$ for every $n\geq1$ (see Proposition~3.1
of~\cite{kn:FIN}).

We claim now that
  \begin{equation}
    \label{eq:sb2}
    \hat Y_n\to Y\,\mbox{ as } n\to\infty
  \end{equation}
almost surely in Skorohod space.

The proof of~(\ref{eq:sb2}) is similar to that of Lemma~\ref{lm:approx}, with
modifications to account for a dependence of $\hat\g_n$ on $n$. (\ref{eq:sb1})
is of course crucial. We indicate the main steps next.

For $n\geq y,\,m\in\N$, let $\d_m$, $\{S_1^m<S_2^m<\ldots\}$ and $L^m_n$ be as in
that proof. We now have that
$\min_{0\leq i\leq L^m_n-1}[\hat\Gamma_n(S^m_{i+1}-)-\hat\Gamma_n(S^m_{i})]>0$
almost surely for $n\geq1$. Define next $\hat\l_n^m:[0,\hat\Gamma_n(S^m_{L^m_n})]\to\R^+$
as follows.
\begin{equation}\label{hln0}
\hat\l_n^m(t)=\frac{\g({\cal Y})}{\hat\g_n({\cal Y}_n)}\,t,\mbox{
if } 0\leq t<\hat\g_n({\cal Y}_n)\,T_0,
\end{equation}
and for $0\leq i\leq L^m_n-1$ and $\hat\Gamma_n(S^m_i)\leq t\leq\hat\Gamma_n(S^m_{i+1})$, let
\begin{equation}
  \label{eq:hln}
  \hat\l_n^m(t)=\begin{cases}\mbox{}\,\,\,\hat\Gamma(S^m_i)\,\,\,\,+\,
  \frac{\hat\Gamma(S^m_{i+1}-)-\hat\Gamma(S^m_{i})}{\hat\Gamma_n(S^m_{i+1}-)
 -\hat\Gamma_n(S^m_{i})}\,\,
  [t-\hat\Gamma_n(S^m_{i})],&
  \mbox{ if } \hat\Gamma_n(S^m_{i})\leq t\leq\hat\Gamma_n(S^m_{i+1}-),\\
  \hat\Gamma(S^m_{i+1}-)
  +\frac{\hat\Gamma(S^m_{i+1})-\hat\Gamma(S^m_{i+1}-)}
  {\hat\Gamma_n(S^m_{i+1})-\hat\Gamma_n(S^m_{i+1}-)}\,
  [t-\hat\Gamma_n(S^m_{i+1}-)],
  &\mbox{ if } \hat\Gamma_n(S^m_{i+1}-)\leq t\leq\hat\Gamma_n(S^m_{i+1}).\end{cases}
\end{equation}
It has the following properties. For all $T>0$, $m\in\N$ and $n\geq m\vee y$
\begin{equation}
  \label{eq:hln1}
  \sup_{0\leq t\leq T}|\hat\l_n^m(t)-t|
  \leq\max\left\{
  |\hat\Gamma(S^m_{i}-)-\hat\Gamma_n(S^m_{i}-)|,\,
  |\hat\Gamma(S^m_{i})-\hat\Gamma_n(S^m_{i})|;\,
  0\leq i\leq L^m_n\right\},
\end{equation}
where $\hat\Gamma_n(0-)\equiv\hat\Gamma(0-)=0$, and the right hand
side of~(\ref{eq:hln1}) vanishes almost surely as $n\to\infty$.
Furthermore,
\begin{equation}
  \label{eq:hln2}
  \sup_{0\leq t\leq
  T}\mbox{dist}\left(\hat Y(\hat\l_n^m(t)),\hat Y_n(t)\right)\leq\d_m,
\end{equation}
since for $t\in[0,T]$, $\hat Y(\hat\l_n^m(t))$ and $\hat Y_n(t)$ coincide
when either one is in $\{1,\ldots,m\}$.

The remainder of the argument follows along the exact same lines as those
in the proof of Lemma~\ref{lm:approx}. \qed

\subsection{Aging}
\label{ssec:age}

Aging results can be viewed as scaling limits for averaged
two-time correlation functions of a given dynamics of a disordered
system. The averaging is with respect to the disorder
distribution. The system should be started at high temperature,
and then abruptly cooled down, evolving thence on at low
temperature. Loosely speaking, aging would amount to the
following. Given a dynamics described by the process $X$ with a
disordered set of parameters $\tau$, the following would be
an aging result.
\begin{equation}
  \label{eq:age}
\lim_{t,t'\to\infty\atop t'/t\to\theta}
\esp\left\{\esp_{\mu}\!\left[\left.\Phi(t,t';X)\right|\tau\right]\right\}
=\RR(\theta),
\end{equation}
where $\mu$ is a measure on state space; $\esp_{\mu}(\cdot|\tau)$ indicates the
expectation with respect to $X$ with initial distribution given by
$\mu$, with parameters fixed at $\tau$; $\Phi(t,t';X)$ is a function of
the piece of trajectory
$X([t,t+t'])=\{X(s),\,s\in[t,t+t']\}$; and $\RR$ is a nontrivial
function of real scaling factor $\theta>0$. The initial
distribution $\mu$ should reflect a high temperature, and the
distribution of the parameters, a low temperature.
See~\cite{kn:BCKM} and references therein.

For a mean field model like the REM-like trap model, there is a volume
dependence, and one must take the infinite volume limit
($n\to\infty$); that should be done before or together with the
time limit. The former is done in~\cite{kn:BD} for
\begin{equation}
  \label{eq:r1}
  \Phi_1(t,t';X)=1\{X(s)=X(t),\,s\in[t,t+t']\}.
\end{equation}
$\mu=\mu_n$ is taken uniformly distributed in $\{1,\ldots,n\}$,
reflecting the high temperature of the initial state, and the tail
parameter $\a<1$ corresponds to the low temperature thence
prevailing.

One other function that is often considered is
\begin{equation}
  \label{eq:r2}
  \Phi_2(t,t';X)=1\{X(t)=X(t+t')\}.
\end{equation}

One could also take the volume and time limits together, using
the scaling limit of Theorem~\ref{teo:scabou}.
For that let us suppose that, for all $t,t'>0$,
$\Phi(t,t';\cdot)$ is almost surely continuous (with respect to
the distribution of $(Y,\g)$). Then, by Theorem~\ref{teo:scabou},
\begin{equation}
  \label{eq:sb3}
\lim_{n\to\infty}
\esp\left\{\esp_{\mu_n}\!\left[\left.\Phi(t,t';\tilde
Y_n)\right|\tau\right]\right\} =
\esp\left\{\esp_{\infty}\!\left[\left.\Phi(t,t';Y)\right|\g\right]\right\},
\end{equation}
with $Y(0)=\infty$, where for $x\in\bN$, $\esp_x(\cdot|\g)$ denotes the
expectation with respect to the distribution of $Y$ started at
$x$, with parameters fixed at $\g$.
$t,t'>0$ are now macroscopic. On those times the dynamics is
already close enough to equilibrium to disallow aging: the right hand side
of~(\ref{eq:sb3}) is {\em not} a function of the ratio $t'/t$
only. To find aging, we should move away from equilibrium, by
taking the further limit as $t,t'\to0$ while $t'/t\to\theta>0$. We
then say that aging takes place in this context (for both the trap
model and the limiting disordered K-process) if
\begin{equation}
  \label{eq:sb4}
\lim_{t,t'\to0\atop t'/t\to\theta}
\esp\left\{\esp_{\infty}\!\left[\left.\Phi(t,t';Y)\right|\g\right]\right\}
=\RR'(\theta)
\end{equation}
exists and $\RR'$ is nontrivial.

\begin{rmk}\label{rmk:tto0}
In taking the volume and time limits as in~\cite{kn:BD}, one
enters what could be termed a (long) microscopic time aging regime
for the trap model, while the latter way of taking those limits
gets one in a (short) macroscopic time aging regime. Our next
result indicates that, at least as far as $\Phi_1$ is concerned,
the two regimes agree.
\end{rmk}

\begin{rmk}\label{rmk:ac}
Instead of scaling time as in~(\ref{eq:sb3}), namely with the
scale of the largest $\tau_x$'s in $K_n$, in view of the further
limit $\lim_{t,t'\to0;\, t'/t\to\theta}$, it is natural to use a
lower order (divergent) scaling. This could be termed  a {\em
mesoscopic} aging regime, and it is the approach of~\cite{kn:AC1}
to establishing aging for the REM-like trap model. As far as 
$\Phi_1$ is concerned, 
the mesoscopic aging regime agrees with the microscopic and
macroscopic regimes;
see~\cite{kn:AC1}.
In the mesoscopic time scaling (as well as in the long microscopic
time scaling), the dynamics itself doesn't have a
limit though.
\end{rmk}

Next we state an aging result for $\Phi$ in a certain class of
functions including the usual examples $\Phi_1$ and $\Phi_2$ and
satisfying some continuity and spatial homogeneity conditions
(which seem natural if one sees this as a mean field model), with
no intention at full generality, however. Let $\Pi$ be the space
of c\`adl\`ag paths on $\bN$, and consider the class of functions
$\Phi:\R^+\times\R^+\times\Pi\to\R$ with the following properties.
\begin{equation}
 \label{eq:sb5}
\Phi(t,t';\zeta)=\Phi(t,t';\zeta([t,t+t'])),
\end{equation}
where $\zeta([t,t+t'])$ is $\zeta$ restricted to $[t,t+t']$,
with the scaling property: for all $t,t'>0$,
\begin{equation}
 \label{eq:sb6}
\Phi(t,t';\zeta)=\Phi(1,t'/t;\zeta^t),
\end{equation}
where $\zeta^t(\cdot)=\zeta(t\cdot)$. Notice that $\Phi_1$ and $\Phi_2$
above have this property. Consider now the following path segments;
for $\theta>0$, $x\in\bN$:
$\eta_{x,\theta}=\eta_{x,\theta}([1,1+\theta))\equiv x$;
$\bar\eta_{x,\theta}=\eta_{x,\theta}([1,1+\theta])\equiv x$. We make
the following further assumptions on $\Phi$:
\begin{equation}
 \label{eq:asf1}
\Phi(1,\theta;\bar\eta_{x,\theta})=\Psi_1(\theta)\,\,\forall x\in\N,
\end{equation}
for some real function $\Psi_1$, that is,
$\Phi(1,\theta;\bar\eta_{x,\theta})$ doesn't depend on $x$ for
finite $x$; for $0<s<\theta$,
\begin{equation}
 \label{eq:asf2}
\Phi(1,\theta;\eta_{x,s}\circ\eta)=\Xi(x,s,\eta)\,\, \forall
x\in\N
\end{equation}
for all segment $\eta=\eta([1+s,1+\theta])$ in $[1+s,1+\theta]$ of
a path in Skorohod space with $\eta(1+s)\ne x$, where  $\circ$
stands for concatenation, and $\Xi$ is a given function with the
following properties. We first give some definitions. For $r>0$,
let $\Pi_u$ be the space of c\`adl\`ag paths in $\bN$ of length
$u$, and for $v>0$ fixed, let ${\cal
X}_v=\cup_{u\in[0,v]}(\{u\}\times\Pi_{v-u})$. Let now $\theta>0$
be fixed. We then have that $\Xi:\bN\times{\cal X}_\theta\to\R$
such that
\begin{eqnarray}\label{eq:asf3a}
&(i)&\Xi \mbox{ is uniformly bounded;}\\
\label{eq:asf3b} &(ii)&\mbox{for all
}\eta\in\Pi_{\theta-s},\,\Xi(x,s,\eta)=\Xi(y,s,\eta)
\mbox{ whenever } x,y\notin\eta;\\
\nn &(iii)&\mbox{for all } x\in\N,\mbox{ the function }\eta\to\Xi(x,s,\eta)
\mbox{ is continuous {\em in the sup norm} on }\\
\label{eq:asf3}
&&\Pi_{\theta-s}\mbox{ for } \eta=\hat\eta_{\infty,\theta-s}\in\Pi_{\theta-s}
\mbox{ with }\hat\eta_{\infty,\theta-s}\equiv\infty, \mbox{ uniformly in }
0<s<\theta.
\end{eqnarray}
\begin{rmk}\label{rmk:f12}
$\Phi_1$ and $\Phi_2$ given in~(\ref{eq:r1},\ref{eq:r2})
satisfy~(\ref{eq:sb5}-\ref{eq:asf3}). Other examples can be
obtained by taking $\Xi(x,s,\eta)=f(s)$ for all $x\in\N$,
$\eta\in\Pi_{\theta-s}$, where $f$ is any continuous
function in $[0,\theta]$.
\end{rmk}
\begin{rmk}\label{rmk:uni}
The uniformity assumption in~(\ref{eq:asf3}) above is for
simplicity. See Remarks~\ref{rmk:cor1} and~\ref{rmk:cor2} below.
\end{rmk}
\begin{rmk}\label{rmk:mf}
The  lack of dependence on finite $x$ assumed in
both~(\ref{eq:asf1}) and~(\ref{eq:asf3b})
is not artificial if one takes into account that the model where
$\Phi$ will be measured is
mean-field, and thus the space coordinate isn't 
relevant. The distinction between finite $x$ and infinity is
nevertheless desirable.
\end{rmk}
For $x\in\N$ and $0<s<\theta$, let
$\Psi_2(s,\theta)=\Xi(x,s,\hat\eta_{\infty,\theta-s})$. Notice that
the latter function doesn't depend on $x\in\N$ by Assumption $(ii)$.
We make the following assumptions on $\Psi_2$, for simplicity: for
all $\theta>0$
\begin{equation}
 \label{eq:asf4}
\Psi_2(0,\theta):=\lim_{s\downarrow0}\Psi_2(s,\theta),\quad
\Psi_2(\theta,\theta):=\lim_{s\uparrow\theta}\Psi_2(s,\theta)\quad\mbox{
exist}
\end{equation}
and
\begin{equation}
 \label{eq:asf5}
\Psi'_2(s,\theta):=\frac{d}{ds}\Psi_2(s,\theta)\in
L_1([0,\theta],\,dx).
\end{equation}

We can now state the main results of this subsection. For simplicity, we make
$t'=\theta t$. We start with a particular case.

\begin{theo}
  \label{teo:age}
For $\g$ as in~(\ref{eq:gf}-\ref{eq:ga}), $\Phi_1$ as
in~(\ref{eq:r1}), and $t,\theta>0$, let
\begin{equation}
  \label{eq:ag0}
\La_t(\theta)=\esp_{\infty}[\Phi_1(t,\theta t;Y)|\g].
\end{equation}
Then almost surely for every $\theta>0$
\begin{equation}
  \label{eq:ag1}
\lim_{t\to0}\La_t(\theta)=\La(\theta),
\end{equation}
where $\La$ is a (nontrivial) function to be exhibited below
(see~(\ref{eq:lam}) and~(\ref{eq:bd})).
\end{theo}

\begin{rmk}\label{rmk:as}
This is an almost sure aging result. The averaged form of~(\ref{eq:sb3})
follows by dominated convergence.
\end{rmk}

\begin{rmk}\label{rmk:bdbf}
As anticipated in Remark~\ref{rmk:tto0} above $\La$ is the same as
the one obtained in~\cite{kn:BD} by taking limits in a different
order and in a different way (see the discussion before and up to
that remark). The computation of the limit in that reference is
not thoroughly rigorous; in~\cite{kn:BF} a rigorous derivation in
the same spirit of~\cite{kn:BD} is performed (in Proposition 2.8
of the former reference).
\end{rmk}

\begin{cor}
  \label{cor:age}
  Let $\Phi$ be as in~(\ref{eq:sb5}-\ref{eq:asf5}).
If~(\ref{eq:ag1}) holds, then
\begin{equation}
  \label{eq:ag2}
\lim_{t\to0} \esp_{\infty}\!\left[\left.\Phi(t,\theta
t;Y)\right|\g\right]
=\Psi_2(0,\theta)+[\Psi_1(\theta)-\Psi_2(\theta,\theta)]\,\La(\theta)
+\int_0^\theta\Psi'_2(s,\theta)\,\La(s)\,ds.
\end{equation}
\end{cor}

\begin{rmk}\label{rmk:cor0}
For the above result, we may take $\g$ fixed such
that~(\ref{eq:ag1}) holds (as well as the assumptions on the
paragraph of~(\ref{eq:ga})).
\end{rmk}

\begin{rmk}\label{rmk:cor1/2}
For both $\Phi_1$ and $\Phi_2$ (see~(\ref{eq:r1},\ref{eq:r2})), we
have $\Psi_1\equiv1$ and $\Psi_2\equiv0$, so, from~(\ref{eq:ag2}),
$\La(\theta)$ is their common aging limit.
\end{rmk}

\begin{rmk}\label{rmk:cor1}
$\La$ turns out to be continuously differentiable in $[0,\infty)$;
we can thus integrate by parts in the right hand side
of~(\ref{eq:ag2}) to obtain that
\begin{equation}
  \label{eq:ag2a}
\lim_{t\to0} \esp_{\infty}\!\left[\left.\Phi(t,\theta
t;Y)\right|\g\right] =\Psi_1(\theta)\,\La(\theta)
-\int_0^\theta\Psi_2(s,\theta)\,\La'(s)\,ds,
\end{equation}
where $\La'(s)=\frac{d}{ds}\La(s)$. For the result in this form we
don't require the uniformity assumption in~(\ref{eq:asf3}), nor
Assumptions (\ref{eq:asf4},\ref{eq:asf5}). See
Remark~\ref{rmk:cor2} below.
\end{rmk}

\begin{rmk}\label{rmk:not}
In the proof of Corollary~\ref{cor:age} below, we will use the
fact that for each $\g$ satisfying the conditions of the paragraph
of~(\ref{eq:sum}) --- in particular, for each $\g$ in a full
measure event, see Remark~\ref{rmk:ga} above
---, and all $t>0$, the distribution of $Y^t$ given $\g$ is the same as that of
$Y$ given $\g^t:=t^{-1}\g$. This follows immediately from the
definition of K-processes (see Definition~\ref{df:kp} and preceding paragraphs).
We thus have that for all such $\g$, and
all bounded measurable function $F$ on Skorohod space,
\begin{equation}
  \label{eq:not}
   \esp\left[\left.F(Y^t)\right|\g\right]
=\esp\left[\left.F(Y)\right|\g^t\right].
\end{equation}
\end{rmk}

\noindent{\bf Proof of Corollary~\ref{cor:age}}  Consider the
conditional expectation on the left hand side of~(\ref{eq:ag2}).
By the scaling property of $\Phi$ (\ref{eq:sb5},\ref{eq:sb6}), we
have that it can be written as
\begin{equation}
  \label{eq:ag3}
\esp\left[\left.\Phi(1,\theta;Y^t)\right|\g\right]
=\esp\left[\left.\Phi(1,\theta;Y)\right|\g^t\right]\, \footnotemark.
\end{equation}
\footnotetext{From now on we write $\P_\infty$ and $\esp_\infty$ as
$\P$ and $\esp$, respectively, using the subscript only for finite
initial points.}
\indent For computing the righ hand side
of~(\ref{eq:ag3}), we first condition on $Y(1)$ and on whether or
not there is a jump of $Y$ in $[1,1+\theta]$, and then if there
is, at which time point it takes place. We get from that and
(\ref{eq:asf1},\ref{eq:asf2})
\begin{eqnarray}\nn&&
\esp\!\left[\left.\Phi(1,\theta;Y)\right|\g^t\right]\\
\label{eq:ag4a}&=& \Psi_1(\theta)\sum_{x\in\N}
\P(Y_1=x|\g^t)\,e^{-\theta
t/\g_x}\\
  \label{eq:ag4b}
&+&\sum_{x\in\N}\P(Y_1=x|\g^t)\int_0^{\theta}\frac{t}{\g_x}\,
e^{-st/\g_x}\,\esp[\Xi(x,s,Y_{[0,\theta-s]})|\g^t]\,ds,
\end{eqnarray}
where the sum can be taken in $\N$ due to Lemma~\ref{lm:zplae},
and we have used time homogeneity of $Y$. \footnote{We have made
notation more compact by substituting parentheses with
subscripts.}

We first note that the sum in~(\ref{eq:ag4a}) equals
$\La_t(\theta)$. Indeed
\begin{eqnarray}\nonumber
\La_t(\theta)
&=&\sum_{x\in\N} \P(Y_t=x|\g)\,\P_x(\mbox{no jump of $Y$ in }[t,t+\theta t]|\g)\\
\label{eq:lt}
&=&\sum_{x\in\N} \P(Y_t=x|\g)\,e^{-\theta
t/\g_x}=\sum_{x\in\N} \P(Y_1=x|\g^t)\,e^{-\theta t/\g_x},
\end{eqnarray}
where we have used the fact alluded to in Remark~\ref{rmk:not}
above in the the third equality. We now write the expression
in~(\ref{eq:ag4b}) as
\begin{equation}
  \label{eq:ag4c1}
\int_0^{\theta}\esp[\Xi(1,s,Y_{[0,\theta-s]})|\g^t]\sum_{x\in\N}
\P(Y_1=x|\g^t)\frac t{\g_x}\,e^{-st/\g_x}\,\,ds
\end{equation}
plus an error that is bounded above by
\begin{equation}
  \label{eq:ag4c2}
\sup_{x\in\N,s\in(0,\theta)}
|\esp[\Xi(1,s,Y_{[0,\theta-s]})|\g^t]-\esp[\Xi(x,s,Y_{[0,\theta-s]})|\g^t]|.
\end{equation}
{}From~(\ref{eq:asf3b}), the absolute value of the difference of
expectations in~(\ref{eq:ag4c2}) can be bounded above by constant
times
\begin{equation}
  \label{eq:ag4c3}
\sup_{x\in\N} \P[1,x\notin Y_{[0,\theta]}|\g^t]=\sup_{x\in\N}
\P[1,x\notin Y_{[0,t\theta]}|\g]
=\sup_{x\in\N}\P[\G^{(1,x)}(t\theta)<\s^{1}_1\vee\s^{x}_1|\g],
\end{equation}
where for $s>0$
\begin{equation}
  \label{eq:ag4c4}
  \G^{(1,x)}(s)=\sum_{y\ne 1,x}\g_x\sum_{i=1}^{N^{(y)}_s}T^{(y)}_i.
\end{equation}
Thus, the right hand side of~(\ref{eq:ag4c3}) is bounded above by
\begin{equation}
  \label{eq:ag4c5}
\P[\G(t\theta)<T'|\g],
\end{equation}
where $T'$ is a continuous random variable independent of $\G$. It
is clear from the fact that $\lim_{s\to0}\G(s)=0$
that~(\ref{eq:ag4c5}) vanishes as $t\to0$.

We thus only have to consider~(\ref{eq:ag4c1}). It can be written
as
\begin{equation}
  \label{eq:ag4c6}
\int_0^{\theta}\Psi_2(s,\theta)\sum_{x\in\N} \P(Y_1=x|\g^t)\frac
t{\g_x}\,e^{-st/\g_x}\,\,ds
=-\int_0^{\theta}\Psi_2(s,\theta)\,\La'_t(s)\,ds,
\end{equation}
where $\La'_t(s)=\frac{d}{ds}\La_t(s)$\footnote{It is a
straightforward exercise to show that the differentiation sign
commutes with the sum.}, plus an error that is bounded above by
\begin{equation}
  \label{eq:ag4c7}
\sup_{s\in(0,\theta)}
|\esp[\Xi(1,s,Y_{[0,\theta-s]})|\g^t]-\Psi_2(s,\theta)|
=\sup_{s\in(0,\theta)}
|\esp[\Xi(1,s,Y_{[0,\theta-s]})|\g^t]-\Xi(1,s,\hat\eta_{\infty,\theta-s})|.
\end{equation}
{}From~(\ref{eq:asf3b},\ref{eq:asf3}), given $\eps>0$, there
exists $\d>0$ such that the difference in~(\ref{eq:ag4c7}) can be
bounded above by constant times
\begin{equation}
  \label{eq:ag4c8}
  \eps+\P\left(\left.\sup_{0\leq s\leq\theta}\mbox{dist}(Y(s),\infty)>\d\right|\g^t\right).
\end{equation}
Since, under $\g^t$, $Y_{[0,\theta]}$ converges in the sup norm to
the identically in $[0,\theta]$ infinity path as $t\to0$, we
conclude that the expression in~(\ref{eq:ag4c7}) vanishes as
$t\to0$.

We are thus left with taking the limit of
\begin{equation}
  \label{eq:ag4d}
\int_0^{\theta}\Psi_2(s,\theta)\,\La'_t(s)\,ds=\Psi_2(\theta,\theta)\,\La_t(\theta)
-\Psi_2(0,\theta)-\int_0^\theta\Psi'_2(s,\theta)\,\La_t(s)\,ds
\end{equation}
as $t\to0$, where we have used the assumptions we made on
$\Psi_2$~(\ref{eq:asf4}-\ref{eq:asf5}) to integrate by parts; note
that $\La_t(0)\equiv1$. Collecting~(\ref{eq:ag4a}-\ref{eq:ag4d})
and the above arguments together with the $L_1$
assumption~(\ref{eq:asf5}) on $\Psi_2'$, the result then follows
by~(\ref{eq:ag1}) and dominated convergence, since $\La_t$ is
bounded (by 1). \qed

\begin{rmk}\label{rmk:cor2}
An alternative, longer argument for the validity of~(\ref{eq:ag2})
in the form~(\ref{eq:ag2a}) {\em for almost every $\g$}, which has
the advantage of requiring neither the uniformity assumption
in~(\ref{eq:asf3}) nor Assumptions (\ref{eq:asf4},\ref{eq:asf5})
--- see Remark~\ref{rmk:cor1} above
--- is to establish the convergence of $\La'_t$ as $t\to0$ to a
(deterministic) function $\La'$ (which turns out to be the
derivative of $\La$) for almost every $\g$. This can be done in an
entirely similar fashion as in the proof of Theorem~\ref{teo:age}
below.  We leave the details for the interested reader.
\end{rmk}

\noindent{\bf Proof of Theorem~\ref{teo:age}}

It is enough to get the result for a fixed $\theta>0$. That we can
find a full measure set of $\g$'s, such that the result holds for
all $\theta>0$ simultaneously, follows from the monotonicity of
$\Phi_1$ and the continuity of $\La$ in $\theta$.

We start with a simpler argument (at this point) for a weaker
result, namely the a.s.~convergence of a (double) Laplace
transform of $\esp[\Phi_1(\cdot,\cdot;Y)|\g]$. This requires the
construction and results of Section~\ref{sec:anal} only. Consider
the function $c_\l(\mu)$ defined in~(\ref{eq:c}). We can represent
it as follows.
\begin{equation}
\label{eq:ag5a}
c_\l(\mu)=
\l\,\mu\int_0^\infty\int_0^\infty e^{-\l s}e^{-\mu t}\esp[\Phi_1(s,t;Y)|\g]\,ds\,dt,
\end{equation}
with $c_\l$ defined in~(\ref{eq:c}). {}For an aging result, it's
natural to take $\mu=\theta\l$, and then take the limit as
$\l\to\infty$. From~(\ref{eq:c1}), we have
\begin{equation}
  \label{eq:ag5b}
c_\l(\l\theta)=\frac{\sum_x\frac{\lambda\g(x)}{1+\lambda\g(x)}
               \frac{\l\theta\g(x)}{1+\l\theta\g(x)}}
                    {\sum_x\frac{\lambda\g(x)}{1+\lambda\g(x)}}.
\end{equation}
Before taking the limit, we note that both sums in~(\ref{eq:ag5b})
can be seen as sums over the increments of the L\'evy process $V$
in $[0,1]$ (see paragraph of~(\ref{eq:tgn}) above) of a function
of the rescaled increments. We thus have by the scale invariance
property of $V$ that the right hand side of~(\ref{eq:ag5b}) for
every $\l>0$ has the same distribution as
\begin{equation}\label{eq:ag5c}
\frac{\sum_{y\in[0,\l^{\a}]}\frac{\g'(y)}{1+\g'(y)}
\frac{\theta\g'(x)}{1+\theta\g'(x)}}
     {\sum_{y\in[0,\l^{\a}]}\frac{\g'(y)}{1+\g'(y)}},
\end{equation}
where the sum is over the increments $\{\g'_x\}$ of $V$ in
$[0,\l^{\a}]$. Now the law of large numbers says that each factor
in the quotient on the right hand side of~(\ref{eq:ag9b3})
converges almost surely to
\begin{equation}\label{eq:ag5d}
 \frac{\esp\sum_{y\in[0,1]}\frac{\g'(y)}{1+\g'(y)}
 \frac{\theta\g'(x)}{1+\theta\g'(x)}}
     {\esp\sum_{y\in[0,1]}\frac{\g'(y)}{1+\g'(y)}}
=\frac{\int_0^\infty\frac{w}{1+w}\,\frac{\theta w}{1+\theta w}\,w^{-1-\a}\,dw}
      {\int_0^\infty\frac{w}{1+w}\,w^{-1-\a}\,dw}
\end{equation}
as $\l\to\infty$.
\begin{rmk}\label{rmk:asconv}
That in principle says that the convergence
of $c_\l(\l\theta)$ as $\l\to\infty$ holds in probability; standard arguments
relying on large deviation estimates for the sums  on the right
hand side of~(\ref{eq:ag5d}) imply convergence almost everywhere.
\end{rmk}

\begin{rmk}\label{rmk:lal}
As we'll see below (in Remark~\ref{rmk:bdbf1}), the expression in
the right hand side of~(\ref{eq:ag5d}) coincides with $\La$ as a
function of $\theta$.
\end{rmk}

We now give a full argument (independent of the above one).
This argument uses the construction and results of
Section~\ref{sec:prob} only.


The argument relies on an estimate for
\begin{equation}
  \label{eq:ag6}
  \P(Y_t=x|\g).
\end{equation}

For $x\in\N$, $\{Y_t=x\}$ can be decomposed in the disjoint union
of
\begin{equation}
  \label{eq:ag7a}
\{\gx(\sx_1)\leq t,\,\gx(\sx_1)+\g_x\,\tix_1>t\},
\end{equation}
where $\gx:=\gx_0$ as in~(\ref{eq:ag7}), and an event where
$\g_x\,\tix_1\leq t$ and $\gx(\sx_2)\leq t$. We thus have
\begin{equation}
  \label{eq:ag8}
  |\P(Y_t=x|\g)-
  \P(\gx(S_1)\leq t,\,\gx(S_1)+\g_x\,\tix_1>t|\g)|
  \leq (1-e^{-\frac t{\g_x}})\,\P(\gx(S_2)\leq t|\g),
\end{equation}
where $S_1$ and $S_2-S_1$ are i.i.d.~rate 1 exponentials which are
independent of all other random variables around.

To establish the result we will prove the two following assertions.
\begin{eqnarray}
\label{eq:ag9a} &\sum_{x\in\N} e^{-\theta t/\g_x}\,\P(\gx(S_1)\leq
t,\,\gx(S_1)+\g_x\,\tix_1>t|\g)\to\La(\theta),&\\
  \label{eq:ag9b}
&\sum_{x\in\N}e^{-\theta t/\g_x}\,\P(\gx(S_2)\leq t|\g)\to0&
\end{eqnarray}
as $t\to0$ for almost every $\g$. We rewrite the probability
in~(\ref{eq:ag9a}) as follows.
\begin{equation}
    \label{eq:ag9c}
\P(\gx(S_1)\leq t|\g)-\P(\gx(S_1)+\g_x\,\tix_1\leq t|\g),
\end{equation}
and note that the second term equals
\begin{equation}
    \label{eq:ag9c1}
\int_0^te^{-(t-s)/\g_x}\,\P(\gx(S_1)\leq
    s|\g)\,ds=\int_0^1e^{-(1-s)t/\g_x}\,\P(\gx(S_1)\leq
    st|\g)\,ds.
\end{equation}
 Substituting in the left hand
side of~(\ref{eq:ag9a}), one sees that in order to prove the
convergence in that display, it is enough to establish
\begin{eqnarray}
\label{eq:ag9d}
&\sum_{x\in\N}e^{-\theta t/\g_x}\,\P(\gx(S_1)\leq t|\g)\to\hat\La(\theta),&\\
  \label{eq:ag9e}
&\int_0^1\sum_{x\in\N}e^{-((1+\theta)-s)t/\g_x}\,\P(\gx(S_1)\leq
    st|\g)\,ds\to\tilde\La(\theta)&
\end{eqnarray}
as $t\to0$ for almost every $\g$, where $\hat\La$ and $\tilde\La$
are functions of $\theta$ only to be given below (see~(\ref{eq:hla})
and~(\ref{eq:tla})); we then have $\La=\hat\La-\tilde\La$.

\begin{rmk}
We note that the left hand sides of~(\ref{eq:ag9d},\ref{eq:ag9e}) are both bounded above
by $\sum_{x\in\N}e^{-\theta t/\g_x}$, which is almost surely finite for every $\theta,t>0$,
since $\sum_{x\in\N}\g_x<\infty$ almost surely. They are thus almost surely finite.
\end{rmk}

We now observe that for almost every $\g$, $\gu\leq\gx\leq\G$ for
all $x\in\N$, where the first domination is a stochastic one
(given $\g$), and follows from the decreasing monotonicity of $\g$.

To get~(\ref{eq:ag9b}), it suffices then to prove that for almost
every $\g$
\begin{equation}\label{eq:ag9b1}
    \P(\gu(S_2)\leq t|\g)\sum_{x\in\N}e^{-\theta t/\g_x}\to0\,\mbox{ as }t\to0.
\end{equation}
For~(\ref{eq:ag9d},\ref{eq:ag9e}), it suffices to prove that for
$i=0,1$
\begin{eqnarray}
\label{eq:ag9d1}
&\P(\gi(S_1)\leq t|\g)\sum_{x\in\N}e^{-\theta t/\g_x}\to\hat\La(\theta),&\\
  \label{eq:ag9e1}
&\int_0^1\P(\gi(S_1)\leq st|\g)
 \sum_{x\in\N}\frac{t}{\g_x}\,e^{-((1+\theta)-s)t/\g_x}\,ds\to\tilde\La(\theta)&
\end{eqnarray}
as $t\to0$, where $\gz=\G$.

The next step is to replace,  for $0<s\leq1$, $i=0,1$, $j=1,2$,
$\P(\gi(S_j)\leq st|\g)$ by constant times $\o_{ij}((st)^{-1})$, where for $r>0$
$\o_{ij}(r):=\esp(\exp\{-r\,\gi(S_j)\}|\g)$. This relies on a
Tauberian theorem (see Theorem 3, Section 5, Chapter XIII
of~\cite{kn:F}), stating that as $t\to0$, the quotient of the
former quantity to the latter one converges to $1/\Ga(j\a)$
provided that for almost every $\g$
\begin{equation}\label{eq:taub}
    \lim_{r\to\infty}\frac{\o_{ij}(qr)}{\o_{ij}(r)}=q^{-j\a}\,
    \mbox{ for all $q>0$,}
\end{equation}
where, for $a>0$, $\Ga(a)=\int_0^\infty t^a e^{-t}\,dt$.

(\ref{eq:taub}) is established in Lemma~\ref{lm:taub} below. From
Lemma~\ref{lm:lap} and~(\ref{eq:la1}), we have that
\begin{equation}\label{eq:ag9b2}
0\leq\o_{12}(t^{-1})\sum_{x\in\N}e^{-\theta t/\g_x}\leq
\frac{\sum_{x}e^{-\theta/t^{-1}\g_x}}
{\left(\sum_{x}\frac{t^{-1}\g_x}{1+t^{-1}\g_x}\right)^2}.
\end{equation}

Arguing as in the sentences above~(\ref{eq:ag5c}), we have that the
right hand side of~(\ref{eq:ag9b2}) for every $t>0$ has the same
distribution as
\begin{equation}\label{eq:ag9b3}
\frac{\sum_{y\in[0,t^{-\a}]}e^{-\theta/\g'_y}}
{\left(\sum_{y\in[0,t^{-\a}]}\frac{\g'_y}{1+\g'_y}\right)^2}=
t^{\a}\,\frac{t^{\a}\sum_{y\in[0,t^{-\a}]}e^{-\theta/\g'_y}}
{\left(t^{\a}\sum_{y\in[0,t^{-\a}]}\frac{\g'_y}{1+\g'_y}\right)^2},
\end{equation}
where the sum is over the increments $\{\g'_x\}$ of $V$ in
$[0,t^{-\a}]$. Now the law of large numbers says that each factor
in the quotient on the right hand side of~(\ref{eq:ag9b3})
converges almost surely to positive finite numbers as $t\to0$. The
extra factor of $t^\a$ in front of that expression then makes it
vanish in that limit. That the same holds for the right hand side
of~(\ref{eq:ag9b2}) follows as in Remark~\ref{rmk:asconv} above.
\medskip

To get~(\ref{eq:ag9d1}), we again need only get the limit for
\begin{equation}\label{eq:ag9d2}
  \o_{i1}(t^{-1})\sum_{x\in\N}e^{-\theta t/\g_x}
\end{equation}
which by Lemma~\ref{lm:lap} and~(\ref{eq:la1}) is bounded above
and below by
\begin{equation}\label{eq:ag9d3}
    \frac{\sum_{x}e^{-\theta/t^{-1}\g_x}}
{k+\sum_{x}\frac{t^{-1}\g_x}{1+t^{-1}\g_x}},
\end{equation}
$k=0$ and $1$, respectively. As in~(\ref{eq:ag9b3}), for every
$t>0$,~(\ref{eq:ag9d3}) has the same distribution as
\begin{equation}\label{eq:ag9d4}
\frac{t^{\a}\sum_{y\in[0,t^{-\a}]}e^{-\theta/\g'_y}}
{kt^{\a}+t^{\a}\sum_{y\in[0,t^{-\a}]}\frac{\g'_y}{1+\g'_y}},
\end{equation}
which by the law of large numbers converges almost surely as
$t\to0$ to
\begin{equation}\label{eq:ag9d5}
\frac{\esp\left(\sum_{y\in[0,1]}e^{-\theta/\g'_y}\right)}
{\esp\left(\sum_{y\in[0,1]}\frac{\g'_y}{1+\g'_y}\right)}=
\frac{\int_0^\infty e^{-\theta/w} w^{-(1+\a)}\,dw}{\int_0^\infty
\frac{w}{1+w} w^{-(1+\a)}\,dw}.
\end{equation}
An analogue of Remark~\ref{rmk:asconv} holds also here.
We thus have from the above that
\begin{equation}\label{eq:hla}
    \hat\La(\theta)=\frac1{\Ga(\a)}\,
    \frac{\int_0^\infty e^{-\theta/w} w^{-(1+\a)}\,dw}{\int_0^\infty
    \frac{w}{1+w} w^{-(1+\a)}\,dw}.
\end{equation}

It remains to get~(\ref{eq:ag9e1}). Since $\P(\gi(S_j)\leq
st|\g)/\o_{ij}((st)^{-1})\to1/\Ga(\a)$ as $t\to0$ uniformly in
$s\in(0,1]$, it suffices to get the limit for
\begin{equation}\label{eq:ag9e2}
\int_0^1\o_{i1}((st)^{-1})\sum_{x\in\N}\frac{t}{\g_x}\,e^{-((1+\theta)-s)t/\g_x}\,ds,
\end{equation}
which by Lemma~\ref{lm:lap} and~(\ref{eq:la1}) reduces to getting
the limit for
\begin{equation}\label{eq:ag9e3}
\int_0^1\frac{\sum_{x}\frac{t}{\g_x}\,e^{-((1+\theta)-s)t/\g_x}}
{k+\sum_{x}\frac{(st)^{-1}\g_x}{1+(st)^{-1}\g_x}}\,ds,
\end{equation}
$k=0,1$. It is clear that the quotient in the above integral, call
it $\hat\La_{s,t}(\theta)$, is bounded above by
$\hat\La_{1,t}(\theta)$, and the latter converges as $t\to0$
almost surely (to $\hat\La(\theta)$, as we just saw). It thus
suffices to establish the almost sure limit of
$\hat\La_{s,t}(\theta)$ as $t\to0$ (independent of $k=0,1$) for every
$s\in(0,1]$. That again works as above: $\hat\La_{s,t}(\theta)$, for every
$t>0$,~(\ref{eq:ag9e3}) has the same distribution as
\begin{equation}\label{eq:ag9e4}
s^\a\,\frac{t^{\a}\sum_{y\in[0,t^{-\a}]}\frac{1}{\g_y'}\,e^{-((1+\theta)-s)/\g'_y}}
{k(st)^{\a}+(st)^{\a}\sum_{y\in[0,(st)^{-\a}]}\frac{\g'_y}{1+\g'_y}}.
\end{equation}
By the law of large numbers, the quotient in~(\ref{eq:ag9e4})
converges almost surely as $t\to0$ to
\begin{equation}\label{eq:ag9e5}
\frac{\esp\left(\sum_{y\in[0,1]}\frac{1}{\g_y'}\,e^{-((1+\theta)-s)/\g'_y}\right)}
{\esp\left(\sum_{y\in[0,1]}\frac{\g'_y}{1+\g'_y}\right)}=
\frac{\int_0^\infty e^{-((1+\theta)-s)/w}
w^{-(2+\a)}\,dw}{\int_0^\infty \frac{w}{1+w} w^{-(1+\a)}\,dw}.
\end{equation}
Again we have an analogue of Remark~\ref{rmk:asconv} here as well.
We thus have from the above that
\begin{equation}\label{eq:tla}
    \tilde\La(\theta)=\frac1{\Ga(\a)}\,
    \frac{\int_0^1\int_0^\infty s^\a\,e^{-((1+\theta)-s)/w} w^{-(2+\a)}\,dw\,ds}{\int_0^\infty
    \frac{w}{1+w} w^{-(1+\a)}\,dw}.
\end{equation}
The result for fixed $\theta>0$ is thus established with
\begin{equation}\label{eq:lam}
    \La=\hat\La-\tilde\La,
\end{equation}
with $\hat\La$, $\tilde\La$ given
in~(\ref{eq:hla}),~(\ref{eq:tla}), respectively.
 \qed

 \begin{rmk}\label{rmk:bdbf1}
It can be shown (e.g., by taking Laplace transforms), that $\La$ thus
obtained coincides with the limit obtained
in~\cite{kn:BD} and~\cite{kn:BF}; in other words
\begin{equation}\label{eq:bd}
    \La(\theta)=\frac{\sin(\pi\a)}{\pi}
    \int_{\frac{\theta}{1+\theta}}^1 s^{-\a}(1-s)^{\a-1}\,ds.
\end{equation}
See Remark~\ref{rmk:bdbf} above. It also coincides with the expression
in the right hand side of~(\ref{eq:ag5d}) as a function of $\theta$.
See Remark~\ref{rmk:lal} above.
 \end{rmk}

 \begin{rmk}\label{rmk:p10bofa}
By~(\ref{eq:lt}), we see that
\begin{equation}\label{eq:bf1}
    \La_t(\theta)=\esp\!\left(\left.e^{-\theta\,t/\g_{\mbox{}_{Y_t}}}\right|\g\right),
\end{equation}
so Theorem~\ref{teo:age} and~(\ref{eq:bd}) establish that, for
almost every $\g$, $t/\g_{\mbox{}_{Y_t}}$ converges in
distribution as $t\to0$ to the random variable $Z$ whose Laplace
transform $\esp(e^{-\theta\,Z})$ is given by the right hand side
of~(\ref{eq:bd}). This is another way of understanding the basic
mechanism for the aging phenomenon in this process ({\em there's
no change for a time of order $t$ when the process has aged $t$
units of time}). It is also a {\em macroscopic} version of the
last assertion of Proposition 2.10 of~\cite{kn:BF}.
\end{rmk}

\begin{rmk}\label{rmk:l11bofa}
Lemma 2.11 of~\cite{kn:BF} establishes the continuity of the
distribution of the random variable $Z$ in
Remark~\ref{rmk:p10bofa}. From~(\ref{eq:hla},\ref{eq:tla}), one
readily finds its density with respect to Lebesgue measure, given
by
\begin{equation}\label{eq:denz}
\frac1{\Ga(\a)\int_0^\infty w^{-\a}\,{(1+w)^{-1}} \,dw}\,\,
z^{\a-1}\int_0^1\a\,s^{\a-1}\,e^{-(1-s)z}ds,\,\,z>0.
\end{equation}
\end{rmk}

\begin{lm}\label{lm:lap} For $i=0,1$, $j=1,2$, and almost every $\g$
\begin{equation}\label{eq:lap}
    \o_{ij}(r)=\left(1+\sum_{x\ne i}\frac{r\g_x}{1+r\g_x}\right)^{-j}.
\end{equation}
\end{lm}

\noindent{\bf Proof} Exercise.

\begin{rmk}\label{rmk:lap}
The condition $x\ne0$ in the sum in~(\ref{eq:lap}) (in the case when $i=0$)
is empty since $x\geq1$. From~(\ref{eq:lap}), we have that for $i=0,1$, $j=1,2$
\begin{equation}\label{eq:la1}
\left(1+\sum_{x}\frac{r\g_x}{1+r\g_x}\right)^{-j}
\leq\o_{ij}(r)\leq
\left(\sum_{x}\frac{r\g_x}{1+r\g_x}\right)^{-j}.
\end{equation}
\end{rmk}

\begin{lm}\label{lm:taub}
(\ref{eq:taub}) holds for almost every $\g$.
\end{lm}

\noindent{\bf Proof} For a fixed $\lambda>0$, and then for all
rational $\lambda>0$, it follows from~(\ref{eq:lap}) and a law of
large number argument as in the proof of Theorem~\ref{teo:age}.
The result for all $\lambda>0$ can be argued from that, using the
monotonicity of $\o_{ij}(\cdot)$ and the continuity of the limit.
\qed

\vspace{.5cm}

\noindent{\bf Acknowledgements} We would like to thank G.~Ben Arous for
letting us have a preliminary version of~\cite{kn:AC1}. L.R.F.~acknowledges
the research grants 307978/2004-4 and 475833/2003-1 by CNPq, and
307978/2004-4, by FAPESP; he also acknowledges support from CNRS. 
The authors benefited from travel grants by the USP-COFECUB and 
the Brazil-France agreements.



\end{document}